\documentclass[11pt,twoside]{article}
\usepackage{amsmath, amsthm, amscd, amsfonts, amssymb, graphicx, color}

\date{}

\begin{document}

\centerline{}

\centerline {\Large{\bf Introduction of retro Banach frame with respect to }}

\centerline {\Large{\bf $b$-linear functional in $n$-Banach space}}

\centerline{}

\newcommand{\mvec}[1]{\mbox{\bfseries\itshape #1}}
\centerline{}
\centerline{\textbf{Prasenjit Ghosh}}
\centerline{Department of Pure Mathematics, University of Calcutta,}
\centerline{35, Ballygunge Circular Road, Kolkata, 700019, West Bengal, India}
\centerline{e-mail: prasenjitpuremath@gmail.com}
\centerline{}
\centerline{\textbf{T. K. Samanta}}
\centerline{Department of Mathematics, Uluberia College,}
\centerline{Uluberia, Howrah, 711315,  West Bengal, India}
\centerline{e-mail: mumpu$_{-}$tapas5@yahoo.co.in}

\newtheorem{Theorem}{\quad Theorem}[section]

\newtheorem{definition}[Theorem]{\quad Definition}

\newtheorem{theorem}[Theorem]{\quad Theorem}

\newtheorem{remark}[Theorem]{\quad Remark}

\newtheorem{corollary}[Theorem]{\quad Corollary}

\newtheorem{note}[Theorem]{\quad Note}

\newtheorem{lemma}[Theorem]{\quad Lemma}

\newtheorem{example}[Theorem]{\quad Example}

\newtheorem{result}[Theorem]{\quad Result}
\newtheorem{conclusion}[Theorem]{\quad Conclusion}

\newtheorem{proposition}[Theorem]{\quad Proposition}

\begin{abstract}
\textbf{\emph{The notion of retro Banach frame with the help of \,$b$-linear functional in \,$n$-Banach spaces is being presented.\,Some properties related to the construction of new retro Banach frame in \,$n$-Banach space have been studied.\,In \,$n$-Banach spaces, some perturbation of retro Banach frames have been discussed.\,Finally, we give a condition in which finite sum of retro Banach frames is a retro Banach frame in \,$n$-Banach space}}
\end{abstract}
{\bf Keywords:}  \emph{Frame, Banach frame, retro Banach frame, $n$-Banach space, \\ \smallskip\hspace{2.2cm}$b$-linear functional.}\\


\section{Introduction}

In Hilbert space, representation of a vector is a most important fact.\,Orthonormal basis is a standard method in Hilbert space to represent a vector and this representation is unique due to the linear independeness of basis.\,A frame is a redundant or linear dependent system for a Hilbert space.\,In mathematical research, frame theory is now a very active field and it has so many applications in wavelet analysis, signal processing, quantum mechanics etc.\,The standard setting for frame is Hilbert space.\,Frame for Hilbert space was introduced by Duffin and Schaeffer\, \cite{Duffin} in 1952.\;Later on, frame theory was popularized by Daubechies et al.\,\cite{Daubechies}.

A sequence \,$\left\{\,f_{\,i}\,\right\}_{i \,=\, 1}^{\infty} \,\subseteq\, H$\, is called a frame for a separable Hilbert space \,$\left(\,H,\, \left<\,\cdot,\, \cdot\,\right>\,\right)$, if there exist positive constants \,$0 \,<\, A \,\leq\, B \,<\, \infty$\, such that
\[ A\; \|\,f\,\|^{\,2} \,\leq\, \sum\limits_{i \,=\, 1}^{\infty}\, \left|\ \left <\,f \,,\, f_{\,i} \, \right >\,\right|^{\,2} \,\leq\, B \,\|\,f\,\|^{\,2}\; \;\text{for all}\; \;f \,\in\, H.\]
A systematic presentation of the frame theory in Hilbert spaces can be found in the book \cite{O}.  

Feichtinger and Groching \cite{F, K} extended the notion of frames to Banach spaces and presented the atomic decomposition for Banach spaces.\,Grochenig also \cite{G} extended Banach frame in more general way in Banach space.\,Thereafter, further development of Banach frame was studied in \cite{H, CH, JK}.\,Retro Banach frame was introduced by P. K. Jain et al.\,\cite{PJ}.\,Stability theorems for Banach frames were studied by Christensen and Heil \cite{CH} and P. K. Jain et al.\,\cite{KJ}  

The notion of linear\;$2$-normed space was introduced by S.\,Gahler \cite{Gahler}.\;The concept of $2$-Banach space is briefly discussed in \cite{White}.\;H.\,Gunawan and Mashadi \cite{Mashadi} developed the generalization of a linear\;$2$-normed space for \,$n \,\geq\, 2$.\,P. Ghosh and T. K. Samanta \cite{PK} introduced the notion of \,$p$-frame relative to \,$b$-linear functional in \,$n$-Banach space.\,They also studied the frames in \,$n$-Hilbert spaces and in their tensor products \cite{Prasenjit, GP} .

In this paper, we introduce the notion of retro Banach frame relative to bounded \,$b$-linear functional in \,$n$-Banach space.\,It has been proved that a \,$n$-Banach space always have a normalized tight retro Banach frame associated to \,$\left(\,a_{\,2},\, \cdots,\, a_{\,n}\,\right)$.\,A necessary and sufficient condition for a sequence in a \,$n$-Banach space to be a retro Banach frame associated to \,$\left(\,a_{\,2},\, \cdots,\, a_{\,n}\,\right)$\, is established.\,Retro Banach frame associated to \,$\left(\,a_{\,2},\, \cdots,\, a_{\,n}\,\right)$\, in Cartesian product of two \,$n$-Banach spaces is also discussed.\,At the end, we describe perturbation results of retro Banach frame associated to \,$(\,a_{\,2},\, \cdots,\, a_{\,n}\,)$.

\section{Preliminaries}

\smallskip\hspace{.6 cm} 

Throughout this paper,\;$E$\, is considered to be a separable Banach space and \,$E^{\,\ast}$, its dual space.\,$\mathcal{B}\,(\,E\,)$\, denotes the space of all bounded linear operators on \,$E$.

\begin{definition}\cite{G}
Let \,$E$\, be a Banach space, \,$E_{d}$\, be a sequence space, which is a Banach space and for which the co-ordinate functionals are continuous.\,Let \,$\left\{\,g_{\,i}\,\right\}_{i \,\in\, I} \,\subseteq\, E^{\,\ast}$\, and \,$S  \,:\, E_{d} \,\to\, E$\, be a bounded linear operator.\,Then the pair \,$\left(\,\{\,g_{\,i}\,\},\,S\,\right)$\, is said to be a Banach frame for \,$E$\, with respect to \,$E_{d}$\, if   
\begin{description}
\item[$(i)$]$\left\{\,g_{\,i}\,(\,f\,)\,\right\} \,\in\, E_{\,d}\; \;\forall\; f \,\in\, E$,
\item[$(ii)$]there exist \,$B \,\geq\, A \,>\, 0$\, such that \,\[A\,\|\,f\,\|_{E} \,\leq\, \left\|\,\left\{\,g_{\,i}\,(\,f\,)\,\right\}\,\right\|_{E_{d}} \,\leq\, B\,\|\,f\,\|_{E}\;\;\forall\, f \,\in\, E,\]
\item[$(iii)$]$S\,\left(\,\left\{\,g_{\,i}\,(\,f\,)\,\right\}\,\right) \,=\, f\; \;\forall\; f \,\in\, E$.
\end{description}
The constants \,$A,\,B$\, are called Banach frame bounds and \,$S$\, is called the reconstruction operator.
\end{definition}

\begin{definition}\cite{PJ}
Let \,$E^{\,\ast}_{\,d}$\, be a Banach space of scalar-valued sequences associated with \,$E^{\,\ast}$\, indexed by \,$\mathbb{N}$.\,Let \,$\left\{\,x_{\,k}\,\right\} \,\subseteq\, E$\, and \,$T \,:\, E^{\,\ast}_{\,d} \,\to\, E^{\,\ast}$\, be given.\,The pair \,$\left(\,\left\{\,x_{\,k}\,\right\},\, T\,\right)$\, is called a retro Banach frame for \,$E^{\,\ast}$\, with respect to \,$E^{\,\ast}_{\,d}$\, if
\begin{description}
\item[$(i)$]$\left\{\,f\,(\,x_{\,k}\,)\,\right\} \,\in\, E^{\,\ast}_{\,d}$, for each \,$f \,\in\, E^{\,\ast}$,
\item[$(ii)$]there exist positive constants \,$A$\, and \,$B$\, with \,$0 \,<\, A \,\leq\, B \,<\, \infty$\, such that \[A\,\|\,f\,\|_{E^{\,\ast}} \,\leq\, \left\|\,\left\{\,f\,(\,x_{\,k}\,)\,\right\}\,\right\|_{E^{\,\ast}_{\,d}} \,\leq\, B\,\|\,f\,\|_{E^{\,\ast}},\; \;f \,\in\, E^{\,\ast},\]
\item[$(iii)$]\,$T$\, is a bounded linear operator such that \,$T\,\left(\,\left\{\,f\,(\,x_{\,k}\,)\,\right\}\,\right) \,=\, f,\; \;f \,\in\, E^{\,\ast}$.
\end{description}
The constants \,$A$\, and \,$B$\, are called frame bounds.\,The operator \,$T$\, is called the reconstruction operator or pre-frame operator.   
\end{definition}

\begin{definition}\cite{Mashadi}
A \,$n$-norm on a linear space \,$X$\, (\,over the field \,$\mathbb{K}$\, of real or complex numbers\,) is a function
\[\left(\,x_{\,1},\, x_{\,2},\, \cdots,\, x_{\,n}\,\right) \,\longmapsto\, \left\|\,x_{\,1},\, x_{\,2},\, \cdots,\, x_{\,n}\,\right\|,\; x_{\,1},\, x_{\,2},\, \cdots,\, x_{\,n} \,\in\, X,\]from \,$X^{\,n}$\, to the set \,$\mathbb{R}$\, of all real numbers such that for every \,$x_{\,1},\, x_{\,2},\, \cdots,\, x_{\,n} \,\in\, X$\, and \,$\alpha \,\in\, \mathbb{K}$,
\begin{description}
\item[$(i)$]\;\; $\left\|\,x_{\,1},\, x_{\,2},\, \cdots,\, x_{\,n}\,\right\| \,=\, 0$\; if and only if \,$x_{\,1},\, \cdots,\, x_{\,n}$\; are linearly dependent,
\item[$(ii)$]\;\;\; $\left\|\,x_{\,1},\, x_{\,2},\, \cdots,\, x_{\,n}\,\right\|$\; is invariant under permutations of \,$x_{\,1},\, x_{\,2},\, \cdots,\, x_{\,n}$,
\item[$(iii)$]\;\;\; $\left\|\,\alpha\,x_{\,1},\, x_{\,2},\, \cdots,\, x_{\,n}\,\right\| \,=\, |\,\alpha\,|\, \left\|\,x_{\,1},\, x_{\,2},\, \cdots,\, x_{\,n}\,\right\|$,
\item[$(iv)$]\;\; $\left\|\,x \,+\, y,\, x_{\,2},\, \cdots,\, x_{\,n}\,\right\| \,\leq\, \left\|\,x,\, x_{\,2},\, \cdots,\, x_{\,n}\,\right\| \,+\,  \left\|\,y,\, x_{\,2},\, \cdots,\, x_{\,n}\,\right\|$.
\end{description}
A linear space \,$X$, together with a n-norm \,$\left\|\,\cdot,\, \cdots,\, \cdot \,\right\|$, is called a linear n-normed space. 
\end{definition}

\begin{definition}\cite{Mashadi}
A sequence \,$\{\,x_{\,k}\,\} \,\subseteq\, X$\, is said to converge to \,$x \,\in\, X$\; if 
\[\lim\limits_{k \to \infty}\,\left\|\,x_{\,k} \,-\, x,\, e_{\,2},\, \cdots,\, e_{\,n} \,\right\| \,=\, 0\; \;\forall\; e_{\,2},\, \cdots,\, e_{\,n} \,\in\, X\]
and it is called a Cauchy sequence if 
\[\lim\limits_{l,\, k \to \infty}\,\left\|\,x_{\,l} \,-\, x_{\,k},\, e_{\,2},\, \cdots,\, e_{\,n}\,\right\| \,=\, 0\; \;\forall\; e_{\,2},\, \cdots,\, e_{\,n} \,\in\, X.\]
The space \,$X$\, is said to be complete or n-Banach space if every Cauchy sequence in this space is convergent in \,$X$.
\end{definition}

\section{Retro Banach frame in $n$-Banach space}

In this section, we first define a bounded\;$b$-linear functional on \,$W \,\times\,\left<\,a_{\,2}\,\right> \,\times\, \cdots \,\times\, \left<\,a_{\,n}\,\right>$, where \,$W$\, be a subspace of \,$X$\, and \,$a_{\,2},\,\cdots,\, a_{\,n} \,\in\, X$\, and then the notion of retro Banach frame in \,$n$-Banach space \,$X$\, is discussed.  

\begin{definition}\label{def1}
Let \,$\left(\,X,\, \|\,\cdot,\, \cdots,\, \cdot\,\| \,\right)$\, be a linear n-normed space and \,$a_{\,2},\, \cdots,\, a_{\,n}$\, be fixed elements in \,$X$.\,Let \,$W$\, be a subspace of \,$X$\, and \,$\left<\,a_{\,i}\,\right>$\, denote the subspaces of \,$X$\, generated by \,$a_{\,i}$, for \,$i \,=\, 2,\, 3,\, \cdots,\,n $.\,Then a map \,$T \,:\, W \,\times\,\left<\,a_{\,2}\,\right> \,\times\, \cdots \,\times\, \left<\,a_{\,n}\,\right> \,\to\, \mathbb{K}$\, is called a b-linear functional defined on \,$W \,\times\, \left<\,a_{\,2}\,\right> \,\times\, \cdots \,\times\, \left<\,a_{\,n}\,\right>$, if for every \,$x,\, y \,\in\, W$\, and \,$k \,\in\, \mathbb{K}$, the following conditions hold:
\begin{description}
\item[$(i)$] $T\,(\,x \,+\, y,\, a_{\,2},\, \cdots,\, a_{\,n}\,) \,=\, T\,(\,x,\, a_{\,2},\, \cdots,\, a_{\,n}\,) \,+\, T\,(\,y,\, a_{\,2},\, \cdots,\, a_{\,n}\,)$
\item[$(ii)$] $T\,(\,k\,x,\, a_{\,2},\, \cdots,\, a_{\,n}\,) \,=\, k\; T\,(\,x,\, a_{\,2},\, \cdots,\, a_{\,n}\,)$. 
\end{description}
A b-linear functional is said to be bounded if there exists a real number \,$M \,>\, 0$\; such that
\[\left|\,T\,(\,x,\, a_{\,2},\, \cdots,\, a_{\,n}\,)\,\right| \,\leq\, M\; \left\|\,x,\, a_{\,2},\, \cdots,\, a_{\,n}\,\right\|\; \;\forall\; x \,\in\, W.\]
The norm of the bounded b-linear functional \,$T$\, is defined by
\[\|\,T\,\| \,=\, \inf\,\left\{\,M \,>\, 0 \,:\, \left|\,T\,(\,x,\, a_{\,2},\, \cdots,\, a_{\,n}\,)\,\right| \,\leq\, M\, \left\|\,x,\, a_{\,2},\, \cdots,\, a_{\,n}\,\right\| \;\forall\; x \,\in\, W\,\right\}.\]
The norm of \,$T$\, can be expressed by any one of the following equivalent formula:
\begin{description}
\item[$(i)$]$\|\,T\,\| \,=\, \sup\,\left\{\,\left|\,T\,(\,x,\, a_{\,2},\, \cdots,\, a_{\,n}\,)\,\right| \;:\; \left\|\,x,\, a_{\,2},\, \cdots,\, a_{\,n}\,\right\| \,\leq\, 1\,\right\}$.
\item[$(ii)$]$\|\,T\,\| \,=\, \sup\,\left\{\,\left|\,T\,(\,x,\, a_{\,2},\, \cdots,\, a_{\,n}\,)\,\right| \;:\; \left\|\,x,\, a_{\,2},\, \cdots,\, a_{\,n}\,\right\| \,=\, 1\,\right\}$.
\item[$(iii)$]$ \|\,T\,\| \,=\, \sup\,\left \{\,\dfrac{\left|\,T\,(\,x,\, a_{\,2},\, \cdots,\, a_{\,n}\,)\,\right|}{\left\|\,x,\, a_{\,2},\, \cdots,\, a_{\,n}\,\right\|} \;:\; \left\|\,x,\, a_{\,2},\, \cdots,\, a_{\,n}\,\right\| \,\neq\, 0\,\right \}$. 
\end{description}
\end{definition}

Some properties of bounded\;$b$-linear functional defined on \,$X \,\times\, \left<\,a_{\,2}\,\right> \,\times\, \cdots \,\times\, \left<\,a_{\,n}\,\right>$\, have been discussed in \cite{TK}.\,For the remaining part of this paper, \,$X$\, denotes the \,$n$-Banach space with respect to the \,$n$-norm \,$\|\,\cdot,\, \cdots,\, \cdot\,\|$\, and \,$X_{F}^{\,\ast}$\, denotes the Banach space of all bounded\;$b$-linear functional defined on \,$X \,\times\, \left<\,a_{\,2}\,\right> \,\times \cdots \,\times\, \left<\,a_{\,n}\,\right>$.

\begin{definition}\label{defn1}
Let \,$X$\, be a \,$n$-Banach space and \,$X^{\,\ast}_{d}$\, be a Banach space of scalar-valued sequences associated with \,$X_{F}^{\,\ast}$\, indexed by \,$\mathbb{N}$.\,Let \,$\left\{\,x_{\,k}\,\right\} \,\subseteq\, X$\, and \,$S \,:\, X^{\,\ast}_{d} \,\to\, X_{F}^{\,\ast}$\, be given.\,Then the pair \,$\left(\,\left\{\,x_{\,k}\,\right\},\, S\,\right)$\, is said to be a retro Banach frame associated to \,$\left(\,a_{\,2},\, \cdots,\, a_{\,n}\,\right)$\, for \,$X_{F}^{\,\ast}$\, with respect to \,$X^{\,\ast}_{d}$\, if 
\begin{description}
\item[$(i)$]$\left\{\,T\,\left(\,x_{\,k},\, a_{\,2},\, \cdots,\, a_{\,n}\,\right)\,\right\} \,\in\, X^{\,\ast}_{d}$, for each \,$T \,\in\, X^{\,\ast}_{F}$.
\item[$(ii)$] There exist constants \,$0 \,<\, A \,\leq\, B \,<\, \infty$\, such that
\begin{equation}\label{eq1}
A\,\|\,T\,\|_{X^{\,\ast}_{F}} \,\leq\, \left\|\,\left\{\,T\,\left(\,x_{\,k},\, a_{\,2},\, \cdots,\, a_{\,n}\,\right)\,\right\}\,\right\|_{X^{\,\ast}_{d}} \,\leq\, B\,\|\,T\,\|_{X^{\,\ast}_{F}}\; \;\forall\; \,T \,\in\, X^{\,\ast}_{F}.
\end{equation}
\item[$(iii)$]$S$\, is a bounded linear operator such that
\[S\,\left(\,\left\{\,T\,\left(\,x_{\,k},\, a_{\,2},\, \cdots,\, a_{\,n}\,\right)\,\right\}\,\right) \,=\, T\; \;\forall\; \,T \,\in\, X^{\,\ast}_{F}.\] 
\end{description}   
\end{definition}

The constants \,$A,\,B$\, are called frame bounds.\,If \,$A \,=\, B$, then \,$\left(\,\left\{\,x_{\,k}\,\right\},\, S\,\right)$\, is called tight retro Banach frame associated to \,$\left(\,a_{\,2},\, \cdots,\, a_{\,n}\,\right)$\, and for \,$A \,=\, B \,=\, 1$, it is called normalized tight retro Banach frame associated to \,$\left(\,a_{\,2},\, \cdots,\, a_{\,n}\,\right)$.\,The inequality (\ref{eq1}), is called the frame inequality for the retro Banach frame associated to \,$\left(\,a_{\,2},\, \cdots,\, a_{\,n}\,\right)$.\,The operator \,$S \,:\, X^{\,\ast}_{d} \,\to\, X_{F}^{\,\ast}$\, is called the reconstruction operator or the pre-frame operator.\,The mapping \,$U \,:\, X^{\,\ast}_{F} \,\to\, X_{d}^{\,\ast}$\, defined by 
\[U\,(\,T\,) \,=\, \left\{\,T\,\left(\,x_{\,k},\, a_{\,2},\, \cdots,\, a_{\,n}\,\right)\,\right\},\; \;T \,\in\, X^{\,\ast}_{F}\]is called the coefficient mapping of the retro Banach frame \,$\left(\,\left\{\,x_{\,k}\,\right\},\, S\,\right)$.\,From inequality (\ref{eq1}), it is easy to verify that \,$U$\, is a bounded linear operator. 

\begin{definition}
A sequence \,$\left\{\,x_{\,k}\,\right\} \,\subseteq\, X$\, is said to be a retro Banach Bessel sequence associated to \,$\left(\,a_{\,2},\, \cdots,\, a_{\,n}\,\right)$\, for \,$X_{F}^{\,\ast}$\, with respect to \,$X^{\,\ast}_{d}$\, if 
\begin{description}
\item[$(i)$]$\left\{\,T\,\left(\,x_{\,k},\, a_{\,2},\, \cdots,\, a_{\,n}\,\right)\,\right\} \,\in\, X^{\,\ast}_{d}$, for each \,$T \,\in\, X^{\,\ast}_{F}$.
\item[$(ii)$]there exists a constant \,$B \,>\, 0$\, such that \[\left\|\,\left\{\,T\,\left(\,x_{\,k},\, a_{\,2},\, \cdots,\, a_{\,n}\,\right)\,\right\}\,\right\|_{X^{\,\ast}_{d}} \,\leq\, B\,\|\,T\,\|_{X^{\,\ast}_{F}}\; \;\forall\; \,T \,\in\, X^{\,\ast}_{F}.\]
\end{description}  
\end{definition}

The constant \,$B$\, is called a retro Banach Bessel bound for the retro Banach Bessel sequence \,$\left\{\,x_{\,k}\,\right\}$\, associated to \,$\left(\,a_{\,2},\, \cdots,\, a_{\,n}\,\right)$.\,Let \,$X_{B}$\, denotes the set of all retro Banach Bessel sequence associated to \,$\left(\,a_{\,2},\, \cdots,\, a_{\,n}\,\right)$\, for \,$X_{F}^{\,\ast}$\, with respect to \,$X^{\,\ast}_{d}$.\,For \,$\left\{\,x_{\,k}\,\right\} \,\in\, X_{B}$, define 
\[\mathcal{R}_{T} \,:\, X_{F}^{\,\ast} \,\to\, X_{d}^{\,\ast}\; \;\text{by}\; \;\mathcal{R}_{T}\,(\,T\,) \,=\, \left\{\,T\,\left(\,x_{\,k},\, a_{\,2},\, \cdots,\, a_{\,n}\,\right)\,\right\},\; \;\text{for all}\; \;T \,\in\, X^{\,\ast}_{F}.\]
Then \,$\mathcal{R}_{T}$\, is a bounded linear operator.\,The operator \,$\mathcal{R}_{T}$\, is called the analysis operator. 

\begin{remark}
Let \,$\left(\,\left\{\,x_{\,k}\,\right\},\, S\,\right)$\, be a retro Banach frame associated to \,$\left(\,a_{\,2},\, \cdots,\, a_{\,n}\,\right)$\, for \,$X_{F}^{\,\ast}$\, with respect to \,$X^{\,\ast}_{d}$\, having frame bounds \,$A$\, and \,$B$.\,Then it is easy to verify that the coefficient mapping \,$U \,:\, X^{\,\ast}_{F} \,\to\, X_{d}^{\,\ast}$\, defined by 
\[U\,(\,T\,) \,=\, \left\{\,T\,\left(\,x_{\,k},\, a_{\,2},\, \cdots,\, a_{\,n}\,\right)\,\right\},\; \;T \,\in\, X^{\,\ast}_{F}\] is a topological isomorphism onto a closed subspace of \,$ X^{\,\ast}_{d}$\, with \,$\|\,U\,\| \,\leq\, B$\, and \,$\left\|\,U^{\,-\, 1}\,\right\| \,\leq\, \dfrac{1}{A}$, where \,$U^{\,-\, 1}$\, exists on the range of \,$U$.
\end{remark}

\begin{proposition}
If for a sequence \,$\left\{\,x_{\,k}\,\right\} \,\subseteq\, X$, the coefficient mapping \,$U$\,  is a topological isomorphism onto \,$X_{d}^{\,\ast}$, then \,$\left(\,\left\{\,x_{\,k}\,\right\},\, S\,\right)$\, is a retro Banach frame associated to \,$\left(\,a_{\,2},\, \cdots,\, a_{\,n}\,\right)$\, for \,$X_{F}^{\,\ast}$\, with respect to \,$X^{\,\ast}_{d}$\, having bounds \,$\left\|\,U^{\,-\, 1}\,\right\|^{\,-\, 1}$\, and \,$\|\,U\,\|$.  
\end{proposition}

\begin{proof}
For each \,$T \,\in\, X^{\,\ast}_{F}$, we have
\[\left\|\,\left\{\,T\,\left(\,x_{\,k},\, a_{\,2},\, \cdots,\, a_{\,n}\,\right)\,\right\}\,\right\|_{X^{\,\ast}_{d}} \,=\, \|\,U\,(\,T\,)\,\|_{X^{\,\ast}_{d}} \,\leq\, \|\,U\,\|\,\|\,T\,\|_{X^{\,\ast}_{F}}.\] 
Also, for \,$T \,\in\, X^{\,\ast}_{F}$, we have 
\begin{align*}
\|\,T\,\|_{X^{\,\ast}_{F}} &\,=\, \left\|\,U^{\,-\, 1}\,\left(\,\left\{\,T\,\left(\,x_{\,k},\, a_{\,2},\, \cdots,\, a_{\,n}\,\right)\,\right\}\,\right)\,\right\|_{X^{\,\ast}_{F}}\\
&\leq\,\|\,U\,\|^{\,-\, 1}\,\left\|\,\left\{\,T\,\left(\,x_{\,k},\, a_{\,2},\, \cdots,\, a_{\,n}\,\right)\,\right\}\,\right\|_{X^{\,\ast}_{d}}. 
\end{align*}
\[\Rightarrow\, \left\|\,U^{\,-\, 1}\,\right\|^{\,-\, 1}\,\|\,T\,\|_{X^{\,\ast}_{F}} \,\leq\, \left\|\,\left\{\,T\,\left(\,x_{\,k},\, a_{\,2},\, \cdots,\, a_{\,n}\,\right)\,\right\}\,\right\|_{X^{\,\ast}_{d}}\; \;\forall\; \,T \,\in\, X^{\,\ast}_{F}.\]
Now, we take \,$S \,=\, U^{\,-\, 1}$.\,Then \,$\left(\,\left\{\,x_{\,k}\,\right\},\, S\,\right)$\, is a retro Banach frame associated to \,$\left(\,a_{\,2},\, \cdots,\, a_{\,n}\,\right)$\, for \,$X_{F}^{\,\ast}$\, with respect to \,$X^{\,\ast}_{d}$\, having frame bounds \,$\left\|\,U^{\,-\, 1}\,\right\|^{\,-\, 1}$\, and \,$\|\,U\,\|$.    
\end{proof}

\begin{definition}
A sequence \,$\left\{\,x_{\,k}\,\right\}$\, in \,$X$\, is said to be total over \,$X^{\,\ast}_{F}$\, if 
\[\left\{\,T \,\in\, X^{\,\ast}_{F} \,:\, T\,\left(\,x_{\,k},\, a_{\,2},\, \cdots,\, a_{\,n}\,\right) \,=\, 0\; \;\forall\; k\,\right\} \,=\, \{\,\theta\,\},\] 
where \,$\theta \,\in\,  X^{\,\ast}_{F}$\, is the null operator.  
\end{definition}

\begin{lemma}\label{lem1}
If a sequence \,$\left\{\,x_{\,k}\,\right\}$\, in \,$X$\, is total over \,$X^{\,\ast}_{F}$, then 
\[X^{\,\ast}_{d_{\,1}} \,=\, \left\{\,\left\{\,T\,\left(\,x_{\,k},\, a_{\,2},\, \cdots,\, a_{\,n}\,\right)\,\right\} \,:\, T \,\in\, X^{\,\ast}_{F}\,\right\}\]is a Banach space with respect to the norm is given by
\begin{equation}\label{eq1.01}
\left\|\,\left\{\,T\,\left(\,x_{\,k},\, a_{\,2},\, \cdots,\, a_{\,n}\,\right)\,\right\}\,\right\|_{X^{\,\ast}_{d_{\,1}}} \,=\, \|\,T\,\|_{X^{\,\ast}_{F}},\; \;T \,\in\, X^{\,\ast}_{F}
\end{equation} 
\end{lemma}

\begin{proof}
It is easy to verify that \,$X^{\,\ast}_{d_{\,1}}$\, is a linear space with respect to the addition and scalar multiplication.\,Also,
\[\left\|\,\left\{\,T\,\left(\,x_{\,k},\, a_{\,2},\, \cdots,\, a_{\,n}\,\right)\,\right\}\,\right\|_{X^{\,\ast}_{d_{\,1}}} \,=\, \|\,T\,\|_{X^{\,\ast}_{F}} \,\geq\, 0\; \;\forall\; \,T \,\in\, X^{\,\ast}_{F}.\]
Now, let \,$\left\|\,\left\{\,T\,\left(\,x_{\,k},\, a_{\,2},\, \cdots,\, a_{\,n}\,\right)\,\right\}\,\right\|_{X^{\,\ast}_{d_{\,1}}} \,=\, 0$.\,Then \,$\|\,T\,\|_{X^{\,\ast}_{F}} \,=\, 0\; \,\Rightarrow\, T \,=\, \theta$.\,This implies that \,$\left\{\,T\,\left(\,x_{\,k},\, a_{\,2},\, \cdots,\, a_{\,n}\,\right)\,\right\} \,=\, 0$.\,On the other hand, if 
\begin{align*}
\left\{\,T\,\left(\,x_{\,k},\, a_{\,2},\, \cdots,\, a_{\,n}\,\right)\,\right\} \,=\, 0 &\,\Rightarrow\, T\,\left(\,x_{\,k},\, a_{\,2},\, \cdots,\, a_{\,n}\,\right) \,=\, 0\; \;\forall\; k.\\
&\Rightarrow\, T \,=\, \theta\; \;[\;\text{since \,$\left\{\,x_{\,k}\,\right\}$ is total over \,$X^{\,\ast}_{F}$}\;]
\end{align*}
This implies that 
\[\left\|\,\left\{\,T\,\left(\,x_{\,k},\, a_{\,2},\, \cdots,\, a_{\,n}\,\right)\,\right\}\,\right\|_{X^{\,\ast}_{d_{\,1}}} \,=\, \|\,T\,\|_{X^{\,\ast}_{F}} \,=\, 0,\; \;T \,\in\, X^{\,\ast}_{F}.\] 
Now, for all \,$\left\{\,T\,\left(\,x_{\,k},\, a_{\,2},\, \cdots,\, a_{\,n}\,\right)\,\right\},\, \left\{\,R\,\left(\,x_{\,k},\, a_{\,2},\, \cdots,\, a_{\,n}\,\right)\,\right\} \,\in\, X^{\,\ast}_{d_{\,1}}$\, and \,$\alpha \,\in\, \mathbb{K}$, we have
\[\left\|\,\alpha\,\left\{\,T\,\left(\,x_{\,k},\, a_{\,2},\, \cdots,\, a_{\,n}\,\right)\,\right\}\,\right\|_{X^{\,\ast}_{d_{\,1}}} \,=\, |\,\alpha\,|\,\left\|\,\left\{\,T\,\left(\,x_{\,k},\, a_{\,2},\, \cdots,\, a_{\,n}\,\right)\,\right\}\,\right\|_{X^{\,\ast}_{d_{\,1}}},\;\text{and}\]
\begin{align*}
&\left\|\,\left\{\,T\,\left(\,x_{\,k},\, a_{\,2},\, \cdots,\, a_{\,n}\,\right)\,\right\} \,+\, \left\{\,R\,\left(\,x_{\,k},\, a_{\,2},\, \cdots,\, a_{\,n}\,\right)\,\right\}\,\right\|_{X^{\,\ast}_{d_{\,1}}}\\
&\leq\, \left\|\,\left\{\,T\,\left(\,x_{\,k},\, a_{\,2},\, \cdots,\, a_{\,n}\,\right)\,\right\}\,\right\|_{X^{\,\ast}_{d_{\,1}}} \,+\, \left\|\,\left\{\,R\,\left(\,x_{\,k},\, a_{\,2},\, \cdots,\, a_{\,n}\,\right)\,\right\}\,\right\|_{X^{\,\ast}_{d_{\,1}}}.
\end{align*}
Thus, \,$X^{\,\ast}_{d_{\,1}}$\, is a normed linear space with respect to the norm given by (\ref{eq1.01}).\\
Let \,$\left\{\,\left\{\,T_{\,i}\,\left(\,x_{\,k},\, a_{\,2},\, \cdots,\, a_{\,n}\,\right)\,\right\}_{\,k}\,\right\}_{\,i}$\, be a Cauchy sequence in \,$X^{\,\ast}_{d_{\,1}}$.\,Then
\begin{align*}
&\lim\limits_{i,\,j \,\to\, \infty}\,\left\|\,\left\{\,T_{\,i}\,\left(\,x_{\,k},\, a_{\,2},\, \cdots,\, a_{\,n}\,\right)\,\right\} \,-\, \left\{\,T_{\,j}\,\left(\,x_{\,k},\, a_{\,2},\, \cdots,\, a_{\,n}\,\right)\,\right\}\,\right\|_{X^{\,\ast}_{d_{\,1}}} \,=\, 0\\
&\Rightarrow\,\lim\limits_{i,\,j \,\to\, \infty}\,\left\|\,\left\{\,\left(\,T_{\,i} \,-\, T_{\,j}\,\right)\,\left(\,x_{\,k},\, a_{\,2},\, \cdots,\, a_{\,n}\,\right)\,\right\}\,\right\|_{X^{\,\ast}_{d_{\,1}}} \,=\, 0.\\
&\Rightarrow\,\lim\limits_{i,\,j \,\to\, \infty}\,\left\|\,T_{\,i} \,-\, T_{\,j}\,\right\|_{X^{\,\ast}_{F}} \,=\, 0\; \;[\;\text{by}\; (\ref{eq1.01})\;]. 
\end{align*} 
This shows that \,$\left\{\,T_{\,i}\,\right\}$\, is a Cauchy sequence in \,$X^{\,\ast}_{F}$.\,Since \,$X^{\,\ast}_{F}$\, is complete, there exists an element \,$T \,\in\, X^{\,\ast}_{F}$\, such that \,$T_{\,i} \,\to\, T$\, as \,$i \,\to\, \infty$.\,Thus, the Cauchy sequence \,$\left\{\,\left\{\,T_{\,i}\,\left(\,x_{\,k},\, a_{\,2},\, \cdots,\, a_{\,n}\,\right)\,\right\}_{\,k}\,\right\}_{\,i}$\, is convergent in \,$X^{\,\ast}_{d_{\,1}}$.\,Hence, \,$X^{\,\ast}_{d_{\,1}}$\, is a Banach space with respect to the norm given by (\ref{eq1.01}).\,This completes the proof. 
\end{proof}

In the next result, we prove that a \,$n$-Banach space always have a normalized tight retro Banach frame associated to \,$\left(\,a_{\,2},\, \cdots,\, a_{\,n}\,\right)$. 

\begin{proposition}
Let \,$X$\, be a \,$n$-Banach space having a retro Banach frame associated to \,$\left(\,a_{\,2},\, \cdots,\, a_{\,n}\,\right)$.\,Then \,$X$\, has a normalized tight retro Banach frame associated to \,$\left(\,a_{\,2},\, \cdots,\, a_{\,n}\,\right)$.  
\end{proposition}

\begin{proof}
Let \,$\left(\,\left\{\,x_{\,k}\,\right\},\, S\,\right)$\, be a retro Banach frame associated to \,$\left(\,a_{\,2},\, \cdots,\, a_{\,n}\,\right)$\, for \,$X_{F}^{\,\ast}$\, with respect to \,$X^{\,\ast}_{d}$.\,Then by the frame inequality (\ref{eq1}), the sequence \,$\left\{\,x_{\,k}\,\right\}$\, is total over \,$X^{\,\ast}_{F}$.\,Therefore, by lemma\,\ref{lem1}, there exists an associated Banach space \,$X_{d_{\,1}}^{\,\ast}$\, with respect to the norm given by (\ref{eq1.01}).\,Define, \,$P \,:\, X^{\,\ast}_{d_{\,1}} \,\to\, X_{F}^{\,\ast}$\, by
\[P\,\left(\,\left\{\,T\,\left(\,x_{\,k},\, a_{\,2},\, \cdots,\, a_{\,n}\,\right)\,\right\}\,\right) \,=\, T,\;  \;T \,\in\, X^{\,\ast}_{F}.\]
Then \,$P$\, is a bounded linear operator such that \,$\left(\,\left\{\,x_{\,k}\,\right\},\, P\,\right)$\, is a normalized tight retro Banach frame associated to \,$\left(\,a_{\,2},\, \cdots,\, a_{\,n}\,\right)$\, for \,$X_{F}^{\,\ast}$\, with respect to \,$X^{\,\ast}_{d}$.   
\end{proof}

\begin{theorem}
Let \,$\left(\,\left\{\,x_{\,k}\,\right\},\, S\,\right)$\, be a retro Banach frame associated to \,$(\,a_{\,2},\, \cdots$, \,$a_{\,n}\,)$\, for \,$X_{F}^{\,\ast}$\, with respect to \,$X^{\,\ast}_{d}$.\,Let \,$\left\{\,y_{\,k}\,\right\} \,\subseteq\, X$\, be such that \,$\left\{\,x_{\,k} \,+\, y_{\,k}\,\right\}$\, is a retro Banach Bessel sequence associated to \,$\left(\,a_{\,2},\, \cdots,\, a_{\,n}\,\right)$\, for \,$X_{F}^{\,\ast}$\, with respect to \,$X^{\,\ast}_{d}$\, having bound \,$K \,<\, \|\,S\,\|^{\,-\, 1}$.\,Then there exists a reconstruction operator \,$P$\, such that \,$\left(\,\left\{\,y_{\,k}\,\right\},\, P\,\right)$\, is a normalized tight retro Banach frame associated to \,$\left(\,a_{\,2},\, \cdots,\, a_{\,n}\,\right)$\, for \,$X_{F}^{\,\ast}$\, with respect to \,$X^{\,\ast}_{d}$.   
\end{theorem}

\begin{proof}
Since \,$\left\{\,x_{\,k} \,+\, y_{\,k}\,\right\}$\, is a retro Banach Bessel sequence associated to \,$(\,a_{\,2},\, \cdots$, \,$a_{\,n}\,)$\, for \,$X_{F}^{\,\ast}$\, with respect to \,$X^{\,\ast}_{d}$\, having bound \,$K$, we have 
\[\left\|\,\left\{\,T\,\left(\,x_{\,k} \,+\, y_{\,k},\, a_{\,2},\, \cdots,\, a_{\,n}\,\right)\,\right\}\,\right\|_{X^{\,\ast}_{d}} \,\leq\, K\,\|\,T\,\|_{X^{\,\ast}_{F}}\; \;\forall\; \,T \,\in\, X^{\,\ast}_{F}.\]
Now, for each \,$T \,\in\, X^{\,\ast}_{F}$, we have
\begin{align*}
&\left(\,\|\,S\,\|^{\,-\, 1} \,-\, K\,\right)\,\|\,T\,\|_{X^{\,\ast}_{F}}\\
& \,\leq\, \left\|\,\left\{\,T\,\left(\,x_{\,k},\, a_{\,2},\, \cdots,\, a_{\,n}\,\right)\,\right\}\,\right\|_{X^{\,\ast}_{d}} \,-\, \left\|\,\left\{\,T\,\left(\,x_{\,k} \,+\, y_{\,k},\, a_{\,2},\, \cdots,\, a_{\,n}\,\right)\,\right\}\,\right\|_{X^{\,\ast}_{d}}\\
&\leq\, \left\|\,\left\{\,T\,\left(\,y_{\,k},\, a_{\,2},\, \cdots,\, a_{\,n}\,\right)\,\right\}\,\right\|_{X^{\,\ast}_{d}}. 
\end{align*}
Thus, the sequence \,$\left\{\,y_{\,k}\,\right\}$\, is total over \,$X^{\,\ast}_{F}$.\,Therefore, by lemma\,\ref{lem1}, there exists an associated Banach space 
\[X^{\,\ast}_{d_{\,1}} \,=\, \left\{\,\left\{\,T\,\left(\,y_{\,k},\, a_{\,2},\, \cdots,\, a_{\,n}\,\right)\,\right\} \,:\, T \,\in\, X^{\,\ast}_{F}\,\right\}\]
equipped with norm given by 
\[\left\|\,\left\{\,T\,\left(\,y_{\,k},\, a_{\,2},\, \cdots,\, a_{\,n}\,\right)\,\right\}\,\right\|_{X^{\,\ast}_{d_{\,1}}} \,=\, \|\,T\,\|_{X^{\,\ast}_{F}},\; \;T \,\in\, X^{\,\ast}_{F}.\]
Define, \,$P \,:\, X^{\,\ast}_{d_{\,1}} \,\to\, X_{F}^{\,\ast}$\, by
\[P\,\left(\,\left\{\,T\,\left(\,y_{\,k},\, a_{\,2},\, \cdots,\, a_{\,n}\,\right)\,\right\}\,\right) \,=\, T,\;  \;T \,\in\, X^{\,\ast}_{F}.\]
Then it is easy to verify that \,$P$\, is a bounded linear operator.\,Hence, \,$\left(\,\left\{\,y_{\,k}\,\right\},\, P\,\right)$\, is a normalized tight retro Banach frame associated to \,$\left(\,a_{\,2},\, \cdots,\, a_{\,n}\,\right)$\, for \,$X_{F}^{\,\ast}$\, with respect to \,$X^{\,\ast}_{d}$.    
\end{proof}

Next, we verify that scalar combinations of two retro Banach frames associated to \,$\left(\,a_{\,2},\, \cdots,\, a_{\,n}\,\right)$\, becomes a retro Banach frame associated to \,$\left(\,a_{\,2},\, \cdots,\, a_{\,n}\,\right)$. 

\begin{theorem}
Let \,$\left(\,\left\{\,x_{\,k}\,\right\},\, S\,\right)$\, and \,$\left(\,\left\{\,y_{\,k}\,\right\},\, S\,\right)$\, be two retro Banach frames associated to \,$\left(\,a_{\,2},\, \cdots,\, a_{\,n}\,\right)$\, for \,$X_{F}^{\,\ast}$\, with respect to \,$X^{\,\ast}_{d}$\, having bounds \,$A,\,B$\, and \,$C,\,D$.\,Then for any scalars \,$\alpha,\,\beta$, \,$\left(\,\left\{\,\alpha\,x_{\,k} \,+\, \beta\,y_{\,k}\,\right\},\, S \,/\, (\,\alpha \,+\, \beta\,)\,\right)$\, is a retro Banach frame associated to \,$\left(\,a_{\,2},\, \cdots,\, a_{\,n}\,\right)$\, for \,$X_{F}^{\,\ast}$\, with respect to \,$X^{\,\ast}_{d}$. 
\end{theorem}

\begin{proof}
For each \,$T \,\in\, X^{\,\ast}_{F}$, we have
\begin{align*}
&\left\|\,\left\{\,T\,\left(\,\alpha\,x_{\,k} \,+\, \beta\,y_{\,k},\, a_{\,2},\, \cdots,\, a_{\,n}\,\right)\,\right\}\,\right\|_{X^{\,\ast}_{d}}\\
&=\,\left\|\,\left\{\,\alpha\,T\,\left(\,x_{\,k},\, a_{\,2},\, \cdots,\, a_{\,n}\,\right) \,+\, \beta\,T\,\left(\,y_{\,k},\, a_{\,2},\, \cdots,\, a_{\,n}\,\right)\,\right\}\,\right\|_{X^{\,\ast}_{d}} \\
&\leq\, |\,\alpha\,|\,\left\|\,\left\{\,T\,\left(\,x_{\,k},\, a_{\,2},\, \cdots,\, a_{\,n}\,\right)\,\right\}\,\right\|_{X^{\,\ast}_{d}} \,+\, |\,\beta\,|\,\left\|\,\left\{\,T\,\left(\,y_{\,k},\, a_{\,2},\, \cdots,\, a_{\,n}\,\right)\,\right\}\,\right\|_{X^{\,\ast}_{d}}\\
&\leq\, \left(\,|\,\alpha\,|\,B \,+\, |\,\beta\,|\,D\,\right)\,\|\,T\,\|_{X^{\,\ast}_{F}}.
\end{align*}
On the other hand,
\begin{align*}
&\left\|\,\left\{\,T\,\left(\,\alpha\,x_{\,k} \,+\, \beta\,y_{\,k},\, a_{\,2},\, \cdots,\, a_{\,n}\,\right)\,\right\}\,\right\|_{X^{\,\ast}_{d}}\\
&\geq\, |\,\alpha\,|\,\left\|\,\left\{\,T\,\left(\,x_{\,k},\, a_{\,2},\, \cdots,\, a_{\,n}\,\right)\,\right\}\,\right\|_{X^{\,\ast}_{d}} \,-\, |\,\beta\,|\,\left\|\,\left\{\,T\,\left(\,y_{\,k},\, a_{\,2},\, \cdots,\, a_{\,n}\,\right)\,\right\}\,\right\|_{X^{\,\ast}_{d}}\\
&\geq\, \left(\,|\,\alpha\,|\,A \,-\, |\,\beta\,|\,C\,\right)\,\|\,T\,\|_{X^{\,\ast}_{F}},\; \;T \,\in\, X^{\,\ast}_{F}.
\end{align*}
Also, for \,$T \,\in\, X^{\,\ast}_{F}$, we have
\[S\,\left(\,\left\{\,T\,\left(\,x_{\,k},\, a_{\,2},\, \cdots,\, a_{\,n}\,\right)\,\right\}\,\right) \,=\, T\; \;\text{and}\; \;S\,\left(\,\left\{\,T\,\left(\,y_{\,k},\, a_{\,2},\, \cdots,\, a_{\,n}\,\right)\,\right\}\,\right) \,=\, T.\] 
Then for \,$T \,\in\, X^{\,\ast}_{F}$, we have
\begin{align*}
&\dfrac{1}{\alpha \,+\, \beta}\,S\,\left(\,\left\{\,T\,\left(\,\alpha\,x_{\,k} \,+\, \beta\,y_{\,k},\, a_{\,2},\, \cdots,\, a_{\,n}\,\right)\,\right\}\,\right)\\
&=\,\dfrac{1}{\alpha \,+\, \beta}\,\left[\,\alpha\,S\,\left(\,\left\{\,T\,\left(\,x_{\,k},\, a_{\,2},\, \cdots,\, a_{\,n}\,\right)\,\right\}\,\right) \,+\, \beta\,S\,\left(\,\left\{\,T\,\left(\,y_{\,k},\, a_{\,2},\, \cdots,\, a_{\,n}\,\right)\,\right\}\,\right)\,\right]\\
&=\,\dfrac{1}{\alpha \,+\, \beta}\,\left(\,\alpha\,T \,+\, \beta\,T\,\right) \,=\, T.  
\end{align*}
Hence, the family \,$\left(\,\left\{\,\alpha\,x_{\,k} \,+\, \beta\,y_{\,k}\,\right\},\, S \,/\, (\,\alpha \,+\, \beta\,)\,\right)$\, is a retro Banach frame associated to \,$(\,a_{\,2}$, \,$\cdots,\, a_{\,n}\,)$\, for \,$X_{F}^{\,\ast}$\, with respect to \,$X^{\,\ast}_{d}$\, having bounds \,$\left(\,|\,\alpha\,|\,A \,-\, |\,\beta\,|\,C\,\right)$\, and \,$\left(\,|\,\alpha\,|\,B \,+\, |\,\beta\,|\,D\,\right)$.   
\end{proof}

In the next theorem, we will see that the sum of two retro Banach frames associated to \,$\left(\,a_{\,2},\, \cdots,\, a_{\,n}\,\right)$\, with different reconstructions operators is also a retro Banach frame associated to \,$\left(\,a_{\,2},\, \cdots,\, a_{\,n}\,\right)$. 

\begin{theorem}
Let \,$\left(\,\left\{\,x_{\,k}\,\right\},\, S\,\right)$\, and \,$\left(\,\left\{\,y_{\,k}\,\right\},\, P\,\right)$\, be two retro Banach frames associated to \,$\left(\,a_{\,2},\, \cdots,\, a_{\,n}\,\right)$\, for \,$X_{F}^{\,\ast}$\, with respect to \,$X^{\,\ast}_{d}$\, having bounds \,$A,\,B$\, and \,$C,\,D$.\,Let \,$R \,:\, X^{\,\ast}_{d} \,\to\, X_{d}^{\,\ast}$\, be a linear homeomorphism such that
\[R\,\left(\,\left\{\,T\,\left(\,x_{\,k},\, a_{\,2},\, \cdots,\, a_{\,n}\,\right)\,\right\}\,\right) \,=\, \left\{\,T\,\left(\,y_{\,k},\, a_{\,2},\, \cdots,\, a_{\,n}\,\right)\,\right\},\; \;T \,\in\, X^{\,\ast}_{F}.\]
Then there exists a reconstruction operator \,$Q \,:\, X^{\,\ast}_{d} \,\to\, X_{F}^{\,\ast}$\, such that the family \,$\left(\,\left\{\,x_{\,k} \,+\, y_{\,k}\,\right\},\,Q\,\right)$\, is a retro Banach frame associated to \,$\left(\,a_{\,2},\, \cdots,\, a_{\,n}\,\right)$\, for \,$X_{F}^{\,\ast}$\, with respect to \,$X^{\,\ast}_{d}$.  
\end{theorem}

\begin{proof}
Let \,$U,\, V$\, be the corresponding coefficient mappings for the retro Banach Bessel sequences \,$\left\{\,x_{\,k}\,\right\}$\, and \,$\left\{\,y_{\,k}\,\right\}$, respectively and \,$I$\, denotes the identity mapping on \,$X^{\,\ast}_{d}$.\,Now, for each \,$T \,\in\, X^{\,\ast}_{F}$, we have
\begin{align*}
&\left\|\,\left\{\,T\,\left(\,x_{\,k} \,+\, y_{\,k},\, a_{\,2},\, \cdots,\, a_{\,n}\,\right)\,\right\}\,\right\|_{X^{\,\ast}_{d}}\\
&=\, \left\|\,\left\{\,T\,\left(\,x_{\,k},\, a_{\,2},\, \cdots,\, a_{\,n}\,\right)\,\right\} \,+\, \left\{\,T\,\left(\,y_{\,k},\, a_{\,2},\, \cdots,\, a_{\,n}\,\right)\,\right\}\,\right\|_{X^{\,\ast}_{d}}\\
&=\, \left\|\,\left\{\,T\,\left(\,x_{\,k},\, a_{\,2},\, \cdots,\, a_{\,n}\,\right)\,\right\} \,+\, R\,\left(\,\left\{\,T\,\left(\,x_{\,k},\, a_{\,2},\, \cdots,\, a_{\,n}\,\right)\,\right\}\,\right)\,\right\|_{X^{\,\ast}_{d}}\\
&\leq\, \|\,I \,+\, R\,\|\,\left\|\,\left\{\,T\,\left(\,x_{\,k},\, a_{\,2},\, \cdots,\, a_{\,n}\,\right)\,\right\}\,\right\|_{X^{\,\ast}_{d}} \,\leq\, B\,\|\,I \,+\, R\,\|\,\|\,T\,\|_{X^{\,\ast}_{F}}.
\end{align*}
Similarly, for each \,$T \,\in\, X^{\,\ast}_{F}$, we have
\[\left\|\,\left\{\,T\,\left(\,x_{\,k} \,+\, y_{\,k},\, a_{\,2},\, \cdots,\, a_{\,n}\,\right)\,\right\}\,\right\|_{X^{\,\ast}_{d}} \,\leq\, D\,\|\,I \,+\, R^{\,-\, 1}\,\|\,\|\,T\,\|_{X^{\,\ast}_{F}}.\]
Thus, for each \,$T \,\in\, X^{\,\ast}_{F}$, we get
\begin{align*}
&\left\|\,\left\{\,T\,\left(\,x_{\,k} \,+\, y_{\,k},\, a_{\,2},\, \cdots,\, a_{\,n}\,\right)\,\right\}\,\right\|_{X^{\,\ast}_{d}}\\
& \,\leq\, \min\,\left\{\,B\,\|\,I \,+\, R\,\|,\, D\,\|\,I \,+\, R^{\,-\, 1}\,\|\,\right\}\,\|\,T\,\|_{X^{\,\ast}_{F}}.
\end{align*}
On the other hand, for each \,$T \,\in\, X^{\,\ast}_{F}$, we have
\begin{align*}
&\left\|\,\left\{\,T\,\left(\,x_{\,k} \,+\, y_{\,k},\, a_{\,2},\, \cdots,\, a_{\,n}\,\right)\,\right\}\,\right\|_{X^{\,\ast}_{d}}\\
&\,\geq\, \left\|\,\left\{\,T\,\left(\,x_{\,k},\, a_{\,2},\, \cdots,\, a_{\,n}\,\right)\,\right\}\,\right\|_{X^{\,\ast}_{d}} \,-\, \left\|\,\left\{\,T\,\left(\,y_{\,k},\, a_{\,2},\, \cdots,\, a_{\,n}\,\right)\,\right\}\,\right\|_{X^{\,\ast}_{d}}\\
&\geq\, A\, \|\,I \,-\, R\,\|\,\|\,T\,\|_{X^{\,\ast}_{F}}.
\end{align*} 
Also, for each \,$T \,\in\, X^{\,\ast}_{F}$, we have  
\[\left\|\,\left\{\,T\,\left(\,x_{\,k} \,+\, y_{\,k},\, a_{\,2},\, \cdots,\, a_{\,n}\,\right)\,\right\}\,\right\|_{X^{\,\ast}_{d}} \,\geq\, C\,\left\|\,R^{\,-\, 1} \,-\, I\,\right\|\,\|\,T\,\|_{X^{\,\ast}_{F}}.\] 
Therefore, for each \,$T \,\in\, X^{\,\ast}_{F}$, we get
\begin{align*}
&\left\|\,\left\{\,T\,\left(\,x_{\,k} \,+\, y_{\,k},\, a_{\,2},\, \cdots,\, a_{\,n}\,\right)\,\right\}\,\right\|_{X^{\,\ast}_{d}}\\
& \,\geq\, \max\left\{\,A\, \|\,I \,-\, R\,\|,\, C\,\left\|\,R^{\,-\, 1} \,-\, I\,\right\|\,\right\}\,\|\,T\,\|_{X^{\,\ast}_{F}}.
\end{align*} 
Now, for \,$T \,\in\, X^{\,\ast}_{F}$, we have   
\begin{align*}
&  R\,\left(\,\left\{\,T\,\left(\,x_{\,k} \,+\, y_{\,k},\, a_{\,2},\, \cdots,\, a_{\,n}\,\right)\,\right\}\,\right)\\
&=\, R\,\left(\,\left\{\,T\,\left(\,x_{\,k},\, a_{\,2},\, \cdots,\, a_{\,n}\,\right)\,\right\}\,\right) \,+\, R\,\left(\,\left\{\,T\,\left(\,y_{\,k},\, a_{\,2},\, \cdots,\, a_{\,n}\,\right)\,\right\}\,\right)\\
&\,=\, (\,I \,+\, R\,)\,\left\{\,T\,\left(\,y_{\,k},\, a_{\,2},\, \cdots,\, a_{\,n}\,\right)\,\right\} \,=\, (\,I \,+\, R\,)\,P^{\,-\, 1}\,T.
\end{align*}
Therefore, if we take \,$Q \,=\, \left(\,(\,I \,+\, R\,)\,P^{\,-\, 1}\,\right)^{\,-\, 1}$, then \,$Q \,:\, X^{\,\ast}_{d} \,\to\, X_{F}^{\,\ast}$\, is a bounded linear operator such that
\[Q\,\left(\,\left\{\,T\,\left(\,x_{\,k} \,+\, y_{\,k},\, a_{\,2},\, \cdots,\, a_{\,n}\,\right)\,\right\}\,\right) \,=\, T\; \;\forall\; \,T \,\in\, X^{\,\ast}_{F}.\]
Hence, \,$\left(\,\left\{\,x_{\,k} \,+\, y_{\,k}\,\right\},\,Q\,\right)$\, is a retro Banach frame associated to \,$\left(\,a_{\,2},\, \cdots,\, a_{\,n}\,\right)$\, for \,$X_{F}^{\,\ast}$\, with respect to \,$X^{\,\ast}_{d}$.       
\end{proof}

The following theorem gives a necessary and sufficient condition for a sequence in a \,$n$-Banach space to be a retro Banach frame associated to \,$\left(\,a_{\,2},\, \cdots,\, a_{\,n}\,\right)$.  

\begin{theorem}
Let \,$\left(\,\left\{\,x_{\,k}\,\right\},\, S\,\right)$\, be a retro Banach frame associated to \,$(\,a_{\,2},\, \cdots$, \,$a_{\,n}\,)$\, for \,$X_{F}^{\,\ast}$\, with respect to \,$X^{\,\ast}_{d}$\, having bounds \,$A,\,B$.\,Let \,$\left\{\,y_{\,k}\,\right\} \,\subseteq\, X$.\,Then there is a reconstruction operator \,$P$\, such that \,$\left(\,\left\{\,y_{\,k}\,\right\},\, P\,\right)$\, is a retro Banach frame associated to \,$\left(\,a_{\,2},\, \cdots,\, a_{\,n}\,\right)$\, for \,$X_{F}^{\,\ast}$\, with respect to \,$X^{\,\ast}_{d}$\, if and only if there exists a constant \,$\lambda \,>\, 0$\, such that
\begin{equation}\label{eq1.1}
\left\|\,R\,\left(\,\left\{\,T\,\left(\,x_{\,k},\, a_{\,2},\, \cdots,\, a_{\,n}\,\right)\,\right\}\,\right)\,\right\|_{X^{\,\ast}_{d}} \,\geq\, \lambda\,\left\|\,\left\{\,T\,\left(\,x_{\,k},\, a_{\,2},\, \cdots,\, a_{\,n}\,\right)\,\right\}\,\right\|_{X^{\,\ast}_{d}},
\end{equation}
for \,$T \,\in\, X^{\,\ast}_{F}$, where \,$R \,:\, X^{\,\ast}_{d} \,\to\, X^{\,\ast}_{d}$\, be a bounded linear operator given by
\[R\,\left(\,\left\{\,T\,\left(\,x_{\,k},\, a_{\,2},\, \cdots,\, a_{\,n}\,\right)\,\right\}\,\right) \,=\, \left\{\,T\,\left(\,y_{\,k},\, a_{\,2},\, \cdots,\, a_{\,n}\,\right)\,\right\},\; \;T \,\in\, X^{\,\ast}_{F}.\]   
\end{theorem}

\begin{proof}
Suppose \,$\left(\,\left\{\,y_{\,k}\,\right\},\, P\,\right)$\, is a retro Banach frame associated to \,$\left(a_{2},\, \cdots,\, a_{n}\right)$\, for \,$X_{F}^{\,\ast}$\, with respect to \,$X^{\,\ast}_{d}$\, having frame bounds \,$C$\, and \,$D$.\,Then
\[C\,\|\,T\,\|_{X^{\,\ast}_{F}} \,\leq\, \left\|\,\left\{\,T\,\left(\,y_{\,k},\, a_{\,2},\, \cdots,\, a_{\,n}\,\right)\,\right\}\,\right\|_{X^{\,\ast}_{d}} \,\leq\, D\,\|\,T\,\|_{X^{\,\ast}_{F}}\; \;\forall\; \,T \,\in\, X^{\,\ast}_{F}.\]
This gives
\begin{align*}
&\left\|\,R\,\left(\,\left\{\,T\,\left(\,x_{\,k},\, a_{\,2},\, \cdots,\, a_{\,n}\,\right)\,\right\}\,\right)\,\right\|_{X^{\,\ast}_{d}} \,\geq\, C\,\|\,T\,\|_{X^{\,\ast}_{F}}\\
&\geq\, \lambda\, \left\|\,\left\{\,T\,\left(\,x_{\,k},\, a_{\,2},\, \cdots,\, a_{\,n}\,\right)\,\right\}\,\right\|_{X^{\,\ast}_{d}},\; \;\text{where \,$\lambda \,=\, \dfrac{C}{B}$}
\end{align*}

Conversely, suppose that (\ref{eq1.1}) holds.\,For each \,$T \,\in\, X^{\,\ast}_{F}$, we have
\begin{align*}
\left\|\,\left\{\,T\,\left(\,y_{\,k},\, a_{\,2},\, \cdots,\, a_{\,n}\,\right)\,\right\}\,\right\|_{X^{\,\ast}_{d}} &\,\geq\, \lambda\, \left\|\,\left\{\,T\,\left(\,x_{\,k},\, a_{\,2},\, \cdots,\, a_{\,n}\,\right)\,\right\}\,\right\|_{X^{\,\ast}_{d}}\\
&\geq\,\lambda\,A\,\|\,T\,\|_{X^{\,\ast}_{F}}.
\end{align*}
On the other hand
\begin{align*}
\left\|\,\left\{\,T\,\left(\,y_{\,k},\, a_{\,2},\, \cdots,\, a_{\,n}\,\right)\,\right\}\,\right\|_{X^{\,\ast}_{d}} & \,=\, \left\|\,R\,\left(\,\left\{\,T\,\left(\,x_{\,k},\, a_{\,2},\, \cdots,\, a_{\,n}\,\right)\,\right\}\,\right)\,\right\|_{X^{\,\ast}_{d}}\\
&\leq\,\|\,R\,\|\,B\,\|\,T\,\|_{X^{\,\ast}_{F}}\; \;\forall\; \,T \,\in\, X^{\,\ast}_{F}. 
\end{align*} 
Let us now define \,$P \,:\, X^{\,\ast}_{d} \,\to\, X_{F}^{\,\ast}$\, by 
\[P\,\left(\,\left\{\,T\,\left(\,y_{\,k},\, a_{\,2},\, \cdots,\, a_{\,n}\,\right)\,\right\}\,\right) \,=\, T,\;  \;T \,\in\, X^{\,\ast}_{F}.\]
Then \,$P$\, is a bounded linear operator such that \,$\left(\,\left\{\,y_{\,k}\,\right\},\, P\,\right)$\, is a retro Banach frame associated to \,$\left(\,a_{\,2},\, \cdots,\, a_{\,n}\,\right)$\, for \,$X_{F}^{\,\ast}$\, with respect to \,$X^{\,\ast}_{d}$.
\end{proof}

We end this section by discussing retro Banach frame associated to \,$\left(\,a_{\,2},\, \cdots,\, a_{\,n}\,\right)$\, in Cartesian product of two \,$n$-Banach spaces.\\

Let \,$\left(\,X,\, \|\,\cdot,\, \cdots,\, \cdot\,\|_{X} \,\right)$\, and \,$\left(\,Y,\, \|\,\cdot,\, \cdots,\, \cdot\,\|_{Y} \,\right)$\, be two \,$n$-Banach spaces.\,Then the Cartesian product of \,$X$\, and \,$Y$\, is denoted by \,$X \,\oplus\, Y$\, and defined to be an \,$n$-Banach space with respect to the \,$n$-norm 
\[\left\|\,x_{\,1} \,\oplus\, y_{\,1},\, x_{\,2} \,\oplus\, y_{\,2},\, \cdots,\, x_{\,n} \,\oplus\, y_{\,n}\,\right\| \,=\, \left\|\,x_{\,1},\, x_{\,2},\, \cdots,\, x_{\,n}\,\right\|_{X} \,+\, \left\|\,y_{\,1},\, y_{\,2},\, \cdots,\, y_{\,n}\,\right\|_{Y},\]
for all \,$x_{\,1} \,\oplus\, y_{\,1},\, x_{\,2} \,\oplus\, y_{\,2},\,\cdots,\, x_{\,n} \,\oplus\, y_{\,n} \,\in\, X \,\oplus\, Y$, and \,$x_{\,1},\, x_{\,2},\, \cdots,\, x_{\,n} \,\in\, X$; \,$\,y_{\,1},\, y_{\,2},\, \cdots,\, y_{\,n} \,\in\, Y$.\,According to the definition (\ref{def1}), consider \,$Y^{\,\ast}_{G}$\, as the Banach space of all bounded\;$b$-linear functional defined on \,$Y \,\times\, \left<\,b_{\,2}\,\right> \,\times \cdots \,\times\, \left<\,b_{\,n}\,\right>$\, and \,$Z_{F \,\oplus\, G}^{\,\ast}$\, as the Banach space of all bounded\;$b$-linear functional defined on \,$X \,\oplus\, Y \,\times\, \left<\,a_{\,2} \,\oplus\, b_{\,2}\,\right> \,\times \cdots \,\times\, \left<\,a_{\,n} \,\oplus\, b_{\,n}\,\right>$, where \,$b_{\,2},\, \cdots,\,b_{\,n} \,\in\, Y$\, and \,$a_{\,2} \,\oplus\, b_{\,2},\, \cdots,\, a_{\,n} \,\oplus\, b_{\,n} \,\in\, X \,\oplus\, Y$\, are fixed elements.\,Now, if \,$T \,\in\, X_{F}^{\,\ast}$\, and \,$U \,\in\, Y^{\,\ast}_{G}$, for all \,$x \,\oplus\, y \,\in\, X \,\oplus\, Y$, we define \,$T \,\oplus\, U \,\in\, Z^{\,\ast}$\, by
\begin{align*}
&(\,T \,\oplus\, U\,)\,\left(\,x \,\oplus\, y,\, a_{\,2} \,\oplus\, b_{\,2},\, \cdots,\, a_{\,n} \,\oplus\, b_{\,n}\,\right)\\
&=\, T\,(\,x,\, a_{\,2},\, \cdots,\, a_{\,n}\,) \,\oplus\, U\,(\,y,\, b_{\,2},\, \cdots,\, b_{\,n}\,)\; \;\forall\; x \,\in\, X,\; y \,\in\, Y. 
\end{align*}  

Let us consider \,$Y^{\,\ast}_{d}$\, \, and \,$Z^{\,\ast}_{d}$\, as the Banach spaces of scalar-valued sequences associated with \,$Y_{G}^{\,\ast}$\, and \,$Z_{F \,\oplus\, G}^{\,\ast}$, respectively.

\begin{theorem}
Let \,$\left(\,\left\{\,x_{\,k}\,\right\},\, S_{X}\,\right)$\, be a retro Banach frame associated to \,$(\,a_{\,2},\, \cdots$, \,$a_{\,n}\,)$\, for \,$X_{F}^{\,\ast}$\, with respect to \,$X^{\,\ast}_{d}$\, having bounds \,$A,\,B$\, and \,$\left(\,\left\{\,y_{\,k}\,\right\},\, S_{Y}\,\right)$\, be a retro Banach frame associated to \,$\left(\,b_{\,2},\, \cdots,\, b_{\,n}\,\right)$\, for \,$Y_{G}^{\,\ast}$\, with respect to \,$Y^{\,\ast}_{d}$\, having bounds \,$C,\,D$.\,Then \,$\left(\,\left\{\,x_{\,k} \,\oplus\, y_{\,k}\,\right\},\, S_{X} \,\oplus\, S_{Y}\,\right)$\, is a retro Banach frame associated to \,$\left(\,a_{\,2} \,\oplus\, b_{\,2},\, \cdots,\, a_{\,n} \,\oplus\, b_{\,n}\,\right)$\, for \,$Z_{F \,\oplus\, G}^{\,\ast}$\, with respect to \,$Z^{\,\ast}_{d}$.
\end{theorem}

\begin{proof}
Since \,$\left(\,\left\{\,x_{\,k}\,\right\},\, S_{X}\,\right)$\, is a retro Banach frame associated to \,$(\,a_{\,2},\, \cdots$, \,$a_{\,n}\,)$\, for \,$X_{F}^{\,\ast}$\, with respect to \,$X^{\,\ast}_{d}$\, and \,$\left(\,\left\{\,y_{\,k}\,\right\},\, S_{Y}\,\right)$\, is a retro Banach frame associated to \,$\left(\,b_{\,2},\, \cdots,\, b_{\,n}\,\right)$\, for \,$Y_{G}^{\,\ast}$\, with respect to \,$Y^{\,\ast}_{d}$, we have
\begin{equation}\label{eq1.2}
A\,\|\,T\,\|_{X^{\,\ast}_{F}} \,\leq\, \left\|\,\left\{\,T\,\left(\,x_{\,k},\, a_{\,2},\, \cdots,\, a_{\,n}\,\right)\,\right\}\,\right\|_{X^{\,\ast}_{d}} \,\leq\, B\,\|\,T\,\|_{X^{\,\ast}_{F}}\; \;\forall\; \,T \,\in\, X^{\,\ast}_{F},
\end{equation} 
\begin{equation}\label{eq1.3}
C\,\|\,R\,\|_{Y^{\,\ast}_{G}} \,\leq\, \left\|\,\left\{\,R\,\left(\,y_{\,k},\, b_{\,2},\, \cdots,\, b_{\,n}\,\right)\,\right\}\,\right\|_{Y^{\,\ast}_{d}} \,\leq\, D\,\|\,R\,\|_{X^{\,\ast}_{F}}\; \;\forall\; \,R \,\in\, Y^{\,\ast}_{G}.
\end{equation}
Adding (\ref{eq1.2}) and (\ref{eq1.3}), we get
\begin{align*}
&A\,\|\,T\,\|_{X^{\,\ast}_{F}} \,+\, C\,\|\,R\,\|_{Y^{\,\ast}_{G}} \\
&\,\leq\, \left\|\,\left\{\,T\,\left(\,x_{\,k},\, a_{\,2},\, \cdots,\, a_{\,n}\,\right)\,\right\}\,\right\|_{X^{\,\ast}_{d}} \,+\, \left\|\,\left\{\,R\,\left(\,y_{\,k},\, b_{\,2},\, \cdots,\, b_{\,n}\,\right)\,\right\}\,\right\|_{Y^{\,\ast}_{d}} \\
&\leq\,B\,\|\,T\,\|_{X^{\,\ast}_{F}} \,+\, D\,\|\,R\,\|_{Y^{\,\ast}_{G}}.\\   
&\Rightarrow\, \min(\,A,\,C\,)\,\left\{\,\|\,T\,\|_{X^{\,\ast}_{F}} \,+\, \,\|\,R\,\|_{Y^{\,\ast}_{G}}\,\right\}\\
& \,\leq\, \left\|\,\left\{\,T\,\left(\,x_{\,k},\, a_{\,2},\, \cdots,\, a_{\,n}\,\right)\,\right\} \,\oplus\, \left\{\,R\,\left(\,y_{\,k},\, b_{\,2},\, \cdots,\, b_{\,n}\,\right)\,\right\}\,\right\|_{Z^{\,\ast}_{d}}\\
& \,\leq\, \max(\,B,\,D\,)\,\left\{\,\|\,T\,\|_{X^{\,\ast}_{F}} \,+\, \,\|\,R\,\|_{Y^{\,\ast}_{G}}\,\right\}.\\
&\Rightarrow\, \min(\,A,\,C\,)\,\left\|\,T \,\oplus\, R\,\right\|_{Z_{F \,\oplus\, G}^{\,\ast}}\\
& \,\leq\, \left\|\,\left\{\,\left(\,T \,\oplus\, R\,\right)\left(\,x_{\,k} \,\oplus\, y_{\,k},\, a_{\,2} \,\oplus\, b_{\,2},\, \cdots,\, a_{\,n} \,\oplus\, b_{\,n}\,\right)\,\right\}\,\right\|_{Z^{\,\ast}_{d}}\\
& \,\leq\, \max(\,B,\,D\,)\,\left\|\,T \,\oplus\, R\,\right\|_{Z_{F \,\oplus\, G}^{\,\ast}}\; \;\forall\; \;T \,\oplus\, R \,\in\, Z_{F \,\oplus\, G}^{\,\ast}.    
\end{align*}
Also, we have
\begin{align*}
&S_{X}\,\left(\,\left\{\,T\,\left(\,x_{\,k},\, a_{\,2},\, \cdots,\, a_{\,n}\,\right)\,\right\}\,\right) \,=\, T,\; \;T \,\in\, X^{\,\ast}_{F},\; \;\text{and}\\
&S_{Y}\,\left(\,\left\{\,R\,\left(\,y_{\,k},\, b_{\,2},\, \cdots,\, b_{\,n}\,\right)\,\right\}\,\right) \,=\, R,\; \;R \,\in\, Y^{\,\ast}_{G} 
\end{align*} 
Now, 
\begin{align*}
&\left(\,S_{X} \,\oplus\, S_{Y}\,\right)\,\left(\,\left\{\,\left(\,T \,\oplus\, R\,\right)\left(\,x_{\,k} \,\oplus\, y_{\,k},\, a_{\,2} \,\oplus\, b_{\,2},\, \cdots,\, a_{\,n} \,\oplus\, b_{\,n}\,\right)\,\right\}\,\right)\\
&\,=\, \left(\,S_{X} \,\oplus\, S_{Y}\,\right)\,\left(\,\left\{\,T\,\left(\,x_{\,k},\, a_{\,2},\, \cdots,\, a_{\,n}\,\right)\,\right\} \,\oplus\, \left\{\,R\,\left(\,y_{\,k},\, b_{\,2},\, \cdots,\, b_{\,n}\,\right)\,\right\}\,\right)\\
&\,=\,S_{X}\,\left(\,\left\{\,T\,\left(\,x_{\,k},\, a_{\,2},\, \cdots,\, a_{\,n}\,\right)\,\right\}\,\right) \,\oplus\, S_{Y}\,\left(\,\left\{\,R\,\left(\,y_{\,k},\, b_{\,2},\, \cdots,\, b_{\,n}\,\right)\,\right\}\,\right)\\
&\,=\, T \,\oplus\, R\; \;\forall\; \;T \,\oplus\, R \,\in\, Z_{F \,\oplus\, G}^{\,\ast}.   
\end{align*}
Hence, the family \,$\left(\,\left\{\,x_{\,k} \,\oplus\, y_{\,k}\,\right\},\, S_{X} \,\oplus\, S_{Y}\,\right)$\, is a retro Banach frame associated to \,$\left(\,a_{\,2} \,\oplus\, b_{\,2},\, \cdots,\, a_{\,n} \,\oplus\, b_{\,n}\,\right)$\, for \,$Z_{F \,\oplus\, G}^{\,\ast}$\, with respect to \,$Z^{\,\ast}_{d}$. 
\end{proof}

\section{Perturbation of retro Banach frame in $n$-Banach space}

In this section, a sufficient condition for the stability of retro Banach frame associated to \,$(\,a_{\,2},\, \cdots,\, a_{\,n}\,)$\, in \,$n$-Banach space under some perturbations is discussed.\,We establish that retro Banach frame associated to \,$(\,a_{\,2},\, \cdots,\, a_{\,n}\,)$\, is stable under perturbation of frame elements by positively confined sequence of scalars.\,Finally, we consider the finite sum of retro Banach frame associated to \,$(\,a_{\,2},\, \cdots,\, a_{\,n}\,)$\, and establish a sufficient condition for the finite sum to be a retro Banach frame associated to \,$(\,a_{\,2},\, \cdots,\, a_{\,n}\,)$\, in \,$n$-Banach space.

Now, we start with a necessary and sufficient condition for the stability of a retro Banach frame associated to \,$\left(\,a_{\,2},\, \cdots,\, a_{\,n}\,\right)$. 

\begin{theorem}\label{thm1}
Let \,$\left(\,\left\{\,x_{\,k}\,\right\},\, S\,\right)$\, be a retro Banach frame associated to \,$(\,a_{\,2},\, \cdots$, \,$a_{\,n}\,)$\, for \,$X_{F}^{\,\ast}$\, with respect to \,$X^{\,\ast}_{d}$\, having bounds \,$A,\,B$.\,Let \,$\left\{\,y_{\,k}\,\right\}$\, be a sequence in \,$X$\, such that $\left\{\,T\,\left(\,y_{\,k},\, a_{\,2},\, \cdots,\, a_{\,n}\,\right)\,\right\} \,\in\, X^{\,\ast}_{d}$, \,$T \,\in\, X^{\,\ast}_{F}$.\,Suppose \,$R \,:\, X_{d}^{\,\ast} \,\to\, X_{d}^{\,\ast}$\, be a bounded linear operator such that
\[R\,\left(\,\left\{\,T\,\left(\,y_{\,k},\, a_{\,2},\, \cdots,\, a_{\,n}\,\right)\,\right\}\,\right) \,=\, \left\{\,T\,\left(\,x_{\,k},\, a_{\,2},\, \cdots,\, a_{\,n}\,\right)\,\right\},\; \;T \,\in\, X_{F}^{\,\ast}.\]
Then there exists a bounded linear operator \,$P \,:\, X^{\,\ast}_{d} \,\to\, X_{F}^{\,\ast}$\, such that \,$\left(\,\left\{\,y_{\,k}\,\right\},\, P\,\right)$\, is a retro Banach frame associated to \,$(\,a_{\,2},\, \cdots,\, a_{\,n}\,)$\, for \,$X_{F}^{\,\ast}$\, with respect to \,$X^{\,\ast}_{d}$\, if and only if there exists a constant \,$K \,>\, 1$\, such that
\begin{align}
&\left\|\,\left\{\,T\,\left(\,x_{\,k} \,-\, y_{\,k},\, a_{\,2},\, \cdots,\, a_{\,n}\,\right)\,\right\}\,\right\|_{X^{\,\ast}_{d}}\nonumber\\
&\leq\, K\,\min\left\{\,\left\|\,\left\{\,T\,\left(\,x_{\,k},\, a_{\,2},\, \cdots,\, a_{\,n}\,\right)\,\right\}\,\right\|_{X^{\,\ast}_{d}},\, \left\|\,\left\{\,T\,\left(\,y_{\,k},\, a_{\,2},\, \cdots,\, a_{\,n}\,\right)\,\right\}\,\right\|_{X^{\,\ast}_{d}}\,\right\}\label{eq1.11} 
\end{align}
\end{theorem}

\begin{proof}
First we suppose that \,$\left(\,\left\{\,y_{\,k}\,\right\},\, P\,\right)$\, is a retro Banach frame associated to \,$(\,a_{\,2},\, \cdots,\, a_{\,n}\,)$.\,Then there exist constants \,$C,\,D \,>\, 0$\, such that
\begin{equation}\label{eq1.21}
A\,\|\,T\,\|_{X^{\,\ast}_{F}} \,\leq\, \left\|\,\left\{\,T\,\left(\,x_{\,k},\, a_{\,2},\, \cdots,\, a_{\,n}\,\right)\,\right\}\,\right\|_{X^{\,\ast}_{d}} \,\leq\, B\,\|\,T\,\|_{X^{\,\ast}_{F}}\; \;\forall\; \,T \,\in\, X^{\,\ast}_{F},
\end{equation}
\begin{equation}\label{eq1.31}
C\,\|\,T\,\|_{X^{\,\ast}_{F}} \,\leq\, \left\|\,\left\{\,T\,\left(\,y_{\,k},\, a_{\,2},\, \cdots,\, a_{\,n}\,\right)\,\right\}\,\right\|_{X^{\,\ast}_{d}} \,\leq\, D\,\|\,T\,\|_{X^{\,\ast}_{F}}\; \;\forall\; \,T \,\in\, X^{\,\ast}_{F}.
\end{equation}
Therefore, for each \,$T \,\in\, X^{\,\ast}_{F}$, we have
\begin{align*}
 & \left\|\,\left\{\,T\,\left(\,x_{\,k} \,-\, y_{\,k},\, a_{\,2},\, \cdots,\, a_{\,n}\,\right)\,\right\}\,\right\|_{X^{\,\ast}_{d}}\\
& \hspace{1cm}\,\leq\, \left\|\,\left\{\,T\,\left(\,x_{\,k},\, a_{\,2},\, \cdots,\, a_{\,n}\,\right)\,\right\}\,\right\|_{X^{\,\ast}_{d}} \,+\, \left\|\,\left\{\,T\,\left(\,y_{\,k},\, a_{\,2},\, \cdots,\, a_{\,n}\,\right)\,\right\}\,\right\|_{X^{\,\ast}_{d}}\\
&\hspace{1cm}\leq\, \left\|\,\left\{\,T\,\left(\,x_{\,k},\, a_{\,2},\, \cdots,\, a_{\,n}\,\right)\,\right\}\,\right\|_{X^{\,\ast}_{d}} \,+\, D\,\|\,T\,\|_{X^{\,\ast}_{F}}\; [\;\text{by (\ref{eq1.31})}\;]\\
&\hspace{1cm} \leq\, \left(\,1 \,+\, \dfrac{D}{A}\,\right)\,\left\|\,\left\{\,T\,\left(\,x_{\,k},\, a_{\,2},\, \cdots,\, a_{\,n}\,\right)\,\right\}\,\right\|_{X^{\,\ast}_{\,d}}\; [\;\text{by (\ref{eq1.21})}\;].  
\end{align*}
Similarly, it can be shown that 
\begin{align*}
&\left\|\,\left\{\,T\,\left(\,x_{\,k} \,-\, y_{\,k},\, a_{\,2},\, \cdots,\, a_{\,n}\,\right)\,\right\}\,\right\|_{X^{\,\ast}_{d}}\\
&\,\leq\, \left(\,1 \,+\, \dfrac{B}{C}\,\right)\,\left\|\,\left\{\,T\,\left(\,y_{\,k},\, a_{\,2},\, \cdots,\, a_{\,n}\,\right)\,\right\}\,\right\|_{X^{\,\ast}_{d}}\; \;\forall\; \,T \,\in\, X^{\,\ast}_{F}.
\end{align*}
Thus, for each \,$T \,\in\, X^{\,\ast}_{F}$, we get
\begin{align*}
&\left\|\,\left\{\,T\,\left(\,x_{\,k} \,-\, y_{\,k},\, a_{\,2},\, \cdots,\, a_{\,n}\,\right)\,\right\}\,\right\|_{X^{\,\ast}_{d}}\\
&\leq\, K\,\min\left\{\,\left\|\,\left\{\,T\,\left(\,x_{\,k},\, a_{\,2},\, \cdots,\, a_{\,n}\,\right)\,\right\}\,\right\|_{X^{\,\ast}_{d}},\, \left\|\,\left\{\,T\,\left(\,y_{\,k},\, a_{\,2},\, \cdots,\, a_{\,n}\,\right)\,\right\}\,\right\|_{X^{\,\ast}_{d}}\,\right\},
\end{align*}
where \,$K \,=\, \max\left\{\,\left(\,1 \,+\, \dfrac{D}{A}\,\right),\, \left(\,1 \,+\, \dfrac{B}{C}\,\right)\,\right\}$.\\

Conversely, suppose that there exists \,$K \,>\, 1$\, such that (\ref{eq1.11}) holds.\,Now, for each \,$T \,\in\, X^{\,\ast}_{F}$, we have
\begin{align*}
&A\,\|\,T\,\|_{X^{\,\ast}_{F}} \,\leq\, \left\|\,\left\{\,T\,\left(\,x_{\,k},\, a_{\,2},\, \cdots,\, a_{\,n}\,\right)\,\right\}\,\right\|_{X^{\,\ast}_{d}}\\
&\leq\, \left\|\,\left\{\,T\,\left(\,x_{\,k} \,-\, y_{\,k},\, a_{\,2},\, \cdots,\, a_{\,n}\,\right)\,\right\}\,\right\|_{X^{\,\ast}_{d}} \,+\, \left\|\,\left\{\,T\,\left(\,y_{\,k},\, a_{\,2},\, \cdots,\, a_{\,n}\,\right)\,\right\}\,\right\|_{X^{\,\ast}_{d}}\\
&\leq\, (\,K \,+\, 1\,)\,\left\|\,\left\{\,T\,\left(\,y_{\,k},\, a_{\,2},\, \cdots,\, a_{\,n}\,\right)\,\right\}\,\right\|_{X^{\,\ast}_{d}}\; [\;\text{by (\ref{eq1.1})}\;]\\
&\Rightarrow\,\dfrac{A}{(\,K \,+\, 1\,)}\,\|\,T\,\|_{X^{\,\ast}_{F}} \,\leq\, \left\|\,\left\{\,T\,\left(\,y_{\,k},\, a_{\,2},\, \cdots,\, a_{\,n}\,\right)\,\right\}\,\right\|_{X^{\,\ast}_{d}}. 
\end{align*}
On the other hand, for each \,$T \,\in\, X^{\,\ast}_{F}$, we have
\begin{align*}
&\left\|\,\left\{\,T\,\left(\,y_{\,k},\, a_{\,2},\, \cdots,\, a_{\,n}\,\right)\,\right\}\,\right\|_{X^{\,\ast}_{d}}\\
& \leq\, \left\|\,\left\{\,T\,\left(\,x_{\,k} \,-\, y_{\,k},\, a_{\,2},\, \cdots,\, a_{\,n}\,\right)\,\right\}\,\right\|_{X^{\,\ast}_{d}} \,+\, \left\|\,\left\{\,T\,\left(\,x_{\,k},\, a_{\,2},\, \cdots,\, a_{\,n}\,\right)\,\right\}\,\right\|_{X^{\,\ast}_{d}}\\
&\leq\, (\,K \,+\, 1\,)\,\left\|\,\left\{\,T\,\left(\,x_{\,k},\, a_{\,2},\, \cdots,\, a_{\,n}\,\right)\,\right\}\,\right\|_{X^{\,\ast}_{d}} \,\leq\, B\,(\,K \,+\, 1\,)\,\|\,T\,\|_{X^{\,\ast}_{F}}.
\end{align*}
Now, take \,$P \,=\, S\,R$.\,Then \,$P \,:\, X^{\,\ast}_{\,d} \,\to\, X_{F}^{\,\ast}$\, is a bounded linear operator such that
\begin{align*}
P\,\left(\,\left\{\,T\,\left(\,y_{\,k},\, a_{\,2},\, \cdots,\, a_{\,n}\,\right)\,\right\}\,\right) &\,=\, S\,R\,\left(\,\left\{\,T\,\left(\,y_{\,k},\, a_{\,2},\, \cdots,\, a_{\,n}\,\right)\,\right\}\,\right)\\
&=\, S\,\left\{\,T\,\left(\,x_{\,k},\, a_{\,2},\, \cdots,\, a_{\,n}\,\right)\,\right\} \,=\, T,\; \;T \,\in\, X^{\,\ast}_{F}. 
\end{align*} 
Thus, \,$\left(\,\left\{\,y_{\,k}\,\right\},\, P\,\right)$\, is a retro Banach frame associated to \,$(\,a_{\,2},\, \cdots,\, a_{\,n}\,)$\, for \,$X_{F}^{\,\ast}$\, with respect to \,$X^{\,\ast}_{d}$.\,This completes the proof. 
\end{proof}

The Theorem \ref{thm1} shows that the stability of retro Banach frame associated to \,$(\,a_{\,2},\, \cdots,\, a_{\,n}\,)$\, depends on the value of \,$K$.\,For large value of \,$K$, the retro Banach frame inequality is lost.\,Therefore, to get optimal frame bounds, we still need stability conditions.\,In the following theorem, we give a sufficient conditions for the stability of a retro Banach frame associated to \,$(\,a_{\,2},\, \cdots,\, a_{\,n}\,)$. 

\begin{theorem}
Let \,$\left(\left\{\,x_{\,k}\,\right\},\, S\,\right)$\, be a retro Banach frame associated to \,$(\,a_{2},\, \cdots,\, a_{n}\,)$\, for \,$X_{F}^{\,\ast}$\, with respect to \,$X^{\,\ast}_{d}$.\,Let \,$\left\{\,y_{\,k}\,\right\} \,\subseteq\, X$\, be such that \,$\left\{\,T\,\left(\,y_{\,k},\, a_{\,2},\, \cdots,\, a_{\,n}\,\right)\,\right\} \,\in\, X^{\,\ast}_{d}$, \,$T \,\in\, X^{\,\ast}_{F}$\, and let \,$U \,:\, X_{F}^{\,\ast} \,\to\, X_{d}^{\,\ast}$\, be the coefficient mapping given by
\[U\,(\,T\,) \,=\, \left\{\,T\,\left(\,x_{\,k},\, a_{\,2},\, \cdots,\, a_{\,n}\,\right)\,\right\},\; \;T \,\in\, X_{F}^{\,\ast}.\]
If there exist positive constants \,$\alpha,\,\beta\;(\,<\, 1\,)$\, and \,$\mu$\, such that 
\begin{description}
\item[$(i)$]\hspace{.2cm}\,$\max\left\{\,\dfrac{\|\,S\,\|\,\left[\,(\,\alpha \,+\, \beta\,)\,\|\,U\,\| \,+\, \mu\,\right]}{(\,1 \,-\, \beta\,)},\, \,\beta\,\right\} \,<\, 1$,
\item[$(ii)$]\hspace{.2cm}$\left\|\,\left\{\,T\,\left(\,x_{\,k} \,-\, y_{\,k},\, a_{\,2},\, \cdots,\, a_{\,n}\,\right)\,\right\}\,\right\|_{X^{\,\ast}_{d}}$
\begin{align*}
&\,\leq\, \alpha\,\left\|\,\left\{\,T\left(\,x_{\,k},\, a_{\,2},\, \cdots,\, a_{\,n}\,\right)\,\right\}\,\right\|_{X^{\,\ast}_{d}} \,+\, \beta\,\left\|\,\left\{\,T\left(\,y_{\,k},\, a_{\,2},\, \cdots,\, a_{\,n}\,\right)\,\right\}\,\right\|_{X^{\,\ast}_{d}}\hspace{7cm}\\
&\hspace{1cm} \,+\, \mu\,\|\,T\,\|_{X^{\,\ast}_{F}},\; \;T \,\in\, X^{\,\ast}_{F},
\end{align*}
\end{description}
then there exists a reconstruction operator \,$P \,:\, X^{\,\ast}_{d} \,\to\, X_{F}^{\,\ast}$\, such that \,$\left(\,\left\{\,y_{\,k}\,\right\},\, P\,\right)$\, is a retro Banach frame associated to \,$(\,a_{\,2},\, \cdots,\, a_{\,n}\,)$\, for \,$X_{F}^{\,\ast}$\, with respect to \,$X^{\,\ast}_{d}$. 
\end{theorem}

\begin{proof}
Let \,$V \,:\, X_{F}^{\,\ast} \,\to\, X_{d}^{\,\ast}$\, be an operator defined by
\[V\,(\,T\,) \,=\, \left\{\,T\,\left(\,y_{\,k},\, a_{\,2},\, \cdots,\, a_{\,n}\,\right)\,\right\},\; \;T \,\in\, X_{F}^{\,\ast}.\]
Using the operators \,$U$\, and \,$V$, condition $(ii)$\, can be written as
\[\left\|\,U\,T \,-\, V\,T\,\right\|_{X^{\,\ast}_{d}} \,\leq\, \alpha\,\left\|\,U\,T\,\right\|_{X^{\,\ast}_{d}} \,+\, \beta\,\left\|\,V\,T\,\right\|_{X^{\,\ast}_{d}} \,+\, \mu\,\|\,T\,\|_{X^{\,\ast}_{F}},\; \;T \,\in\, X_{F}^{\,\ast}.\]
Thus, for \,$T \,\in\, X_{F}^{\,\ast}$, we have
\begin{align*}
\left\|\,\left\{\,T\left(\,y_{\,k},\, a_{2},\, \cdots,\, a_{n}\,\right)\,\right\}\,\right\|_{X^{\,\ast}_{d}} \,=\, \left\|\,V\,T\,\right\|_{X^{\,\ast}_{d}} &\,\leq\, \left\|\,U\,T \,-\, V\,T\,\right\|_{X^{\,\ast}_{d}} \,+\, \left\|\,U\,T\,\right\|_{X^{\,\ast}_{d}}\\
&\,\leq\, \dfrac{(\,1 \,+\, \alpha\,)\,\|\,U\,\| \,+\, \mu}{1 \,-\, \beta}\,\|\,T\,\|_{X^{\,\ast}_{F}}. 
\end{align*}
Therefore, \,$V$\, is a bounded linear operator such that
\[\left\|\,U\,T \,-\, V\,T\,\right\|_{X^{\,\ast}_{d}} \,\leq\, \dfrac{(\,\alpha \,+\, \beta\,)\,\|\,U\,\| \,+\, \mu}{1 \,-\, \beta}\,\|\,T\,\|_{X^{\,\ast}_{F}}\; \;\forall\; T \,\in\, X_{F}^{\,\ast}.\] 
Now, 
\[\left\|\,I_{F} \,-\, S\,V\,\right\| \,\leq\, \|\,S\,\|\,\|\,U \,-\, V\,\| \,\leq\, \dfrac{\|\,S\,\|\,\left[\,(\,\alpha \,+\, \beta\,)\,\|\,U\,\| \,+\, \mu\,\right]}{1 \,-\, \beta} \,<\, 1.\]
This shows that \,$S\,V$\, is an invertible operator with satisfying 
\[\|\,(\,S\,V\,)^{\,-\, 1}\,\| \,\leq\, \dfrac{1}{1 \,-\, \dfrac{\|\,S\,\|\,\left[\,(\,\alpha \,+\, \beta\,)\,\|\,U\,\| \,+\, \mu\,\right]}{1 \,-\, \beta}}.\]
Now, take \,$P \,=\, (\,S\,V\,)^{\,-\, 1}\,S$.\,Then \,$P\,V \,=\, I_{F}$, where \,$I_{F}$\, is the identity operator on \,$X_{F}^{\,\ast}$.\,Thus, \,$P \,:\, X^{\,\ast}_{\,d} \,\to\, X_{F}^{\,\ast}$\, is a bounded linear operator such that
\[P\,\left(\,\left\{\,T\,\left(\,y_{\,k},\, a_{\,2},\, \cdots,\, a_{\,n}\,\right)\,\right\}\,\right) \,=\, T,\; \;T \,\in\, X^{\,\ast}_{F}.\]
Now, for \,$T \,\in\, X_{F}^{\,\ast}$, we have
\begin{align*}
\|\,T\,\|_{X^{\,\ast}_{F}} \,=\, \|\,P\,V\,T\,\|_{X^{\,\ast}_{F}} &\,\leq\, \|\,P\,\|\,\|\,V\,T\,\|_{X^{\,\ast}_{d}}\\
& \,\leq\, \dfrac{\|\,S\,\|}{1 \,-\, \dfrac{\|\,S\,\|\,\left[\,(\,\alpha \,+\, \beta\,)\,\|\,U\,\| \,+\, \mu\,\right]}{1 \,-\, \beta}}\,\|\,V\,T\,\|_{X^{\,\ast}_{d}}.
\end{align*}
This implies that
\begin{align*}
& \|\,S\,\|^{\,-\, 1}\,\left(\,1 \,-\, \dfrac{\left[\,(\,\alpha \,+\, \beta\,)\,\|\,U\,\| \,+\, \mu\,\right]\,\|\,S\,\|}{1 \,-\, \beta}\,\right)\,\|\,T\,\|_{X^{\,\ast}_{F}}\\
&  \leq\, \left\|\,\left\{\,T\left(\,y_{\,k},\, a_{\,2},\, \cdots,\, a_{\,n}\,\right)\right\}\,\right\|_{X^{\,\ast}_{d}}\; \;\forall\; T \,\in\, X_{F}^{\,\ast}.
\end{align*}
Hence, \,$\left(\,\left\{\,y_{\,k}\,\right\},\, P\,\right)$\, is a retro Banach frame associated to \,$(\,a_{\,2},\, \cdots,\, a_{\,n}\,)$\, for \,$X_{F}^{\,\ast}$\, with respect to \,$X^{\,\ast}_{d}$.\,This completes the proof.    
\end{proof}

Next, we give a stability condition of a retro Banach frame associated to \,$(\,a_{\,2},\, \cdots$, \,$a_{\,n}\,)$\, by uniformly scaled version of a given retro Banach Bessel sequence associated to \,$(\,a_{\,2},\, \cdots,\, a_{\,n}\,)$. 

\begin{theorem}\label{th1}
Let \,$\left(\left\{\,x_{\,k}\,\right\},\, S\,\right)$\, be a retro Banach frame associated to \,$(\,a_{2},\, \cdots,\, a_{n}\,)$\, for \,$X_{F}^{\,\ast}$\, with respect to \,$X^{\,\ast}_{d}$\, having bounds \,$A,\,B$.\,Let \,$\left\{\,y_{\,k}\,\right\}$\, be a sequence in \,$X$\, such that \,$\left\{\,T\,\left(\,y_{\,k},\, a_{\,2},\, \cdots,\, a_{\,n}\,\right)\,\right\} \,\in\, X^{\,\ast}_{d}$, \,$T \,\in\, X^{\,\ast}_{F}$\, and for some constant \,$K \,>\, 0$
\[ \left\|\,\left\{\,T\,\left(\,y_{\,k},\, a_{\,2},\, \cdots,\, a_{\,n}\,\right)\,\right\}\,\right\|_{X^{\,\ast}_{d}} \,\leq\, K\,\|\,T\,\|_{X^{\,\ast}_{F}},\; \;\text{for all}\; \;T \,\in\, X^{\,\ast}_{F}.\]
Then for any non-zero constant \,$\lambda$\, with \,$|\,\lambda\,| \,<\, \dfrac{\|\,S\,\|^{\,-\, 1}}{K}$, there exists a reconstruction operator \,$P \,:\, X^{\,\ast}_{d} \,\to\, X_{F}^{\,\ast}$\, such that \,$\left(\,\left\{\,x_{\,k} \,+\, \lambda\,y_{\,k}\,\right\},\, P\,\right)$\, is a retro Banach frame associated to \,$(\,a_{\,2},\, \cdots,\, a_{\,n}\,)$\, for \,$X_{F}^{\,\ast}$\, with respect to \,$X^{\,\ast}_{d}$\, having frame bounds \,$\left(\,\|\,S\,\|^{\,-\, 1} \,-\, |\,\lambda\,|\,K\,\right)$\, and \,$\left(\,B \,+\, |\,\lambda\,|\,K\,\right)$.
\end{theorem}

\begin{proof}
Let \,$V \,:\, X_{F}^{\,\ast} \,\to\, X_{d}^{\,\ast}$\, be a bounded linear operator defined by
\[V\,(\,T\,) \,=\, \left\{\,T\,\left(\,y_{\,k},\, a_{\,2},\, \cdots,\, a_{\,n}\,\right)\,\right\},\; \;T \,\in\, X_{F}^{\,\ast},\] 
and \,$U \,:\, X_{F}^{\,\ast} \,\to\, X_{d}^{\,\ast}$\, be a bounded linear operator given by
\[U\,(\,T\,) \,=\, \left\{\,T\,\left(\,x_{\,k},\, a_{\,2},\, \cdots,\, a_{\,n}\,\right)\,\right\},\; \;T \,\in\, X_{F}^{\,\ast}.\]  
Then it is easy to verify that \,$\left\{\,T\,\left(\,x_{\,k} \,+\, \lambda\,y_{\,k},\, a_{\,2},\, \cdots,\, a_{\,n}\,\right)\,\right\} \,\in\, X^{\,\ast}_{d}$, for all \,$T \,\in\, X^{\,\ast}_{F}$.\,Now, for each \,$T \,\in\, X^{\,\ast}_{F}$, we have
\begin{align*}
\left\|\,U\,T \,+\, \lambda\,V\,T\,\right\|_{X^{\,\ast}_{d}} & \,=\, \left\|\,\left\{\,T\,\left(\,x_{\,k} \,+\, \lambda\,y_{\,k},\, a_{\,2},\, \cdots,\, a_{\,n}\,\right)\,\right\}\,\right\|_{X^{\,\ast}_{d}}\\
& \,\leq\, \left\|\,\left\{\,T\,\left(\,x_{\,k},\, a_{\,2},\, \cdots,\, a_{\,n}\,\right)\,\right\}\,\right\|_{X^{\,\ast}_{d}}\\
&\hspace{1cm} \,+\, |\,\lambda\,|\,\left\|\,\left\{\,T\,\left(\,y_{\,k},\, a_{\,2},\, \cdots,\, a_{\,n}\,\right)\,\right\}\,\right\|_{X^{\,\ast}_{d}}\\
& \leq\, \left(\,B \,+\, |\,\lambda\,|\,K\,\right)\,\|\,T\,\|_{X^{\,\ast}_{F}}.
\end{align*}
On the other hand, for each \,$T \,\in\, X^{\,\ast}_{F}$, we have
\begin{align*}
\left(\,\|\,S\,\|^{\,-\, 1} \,-\, |\,\lambda\,|\,K\,\right)\,\|\,T\,\|_{X^{\,\ast}_{F}} & \,\leq\, \left\|\,\left\{\,T\,\left(\,x_{\,k},\, a_{\,2},\, \cdots,\, a_{\,n}\,\right)\,\right\}\,\right\|_{X^{\,\ast}_{d}}\\
&\hspace{1cm} \,-\, |\,\lambda\,|\,\left\|\,\left\{\,T\,\left(\,y_{\,k},\, a_{\,2},\, \cdots,\, a_{\,n}\,\right)\,\right\}\,\right\|_{X^{\,\ast}_{d}}\\
&\leq\, \left\|\,\left\{\,T\,\left(\,x_{\,k} \,+\, \lambda\,y_{\,k},\, a_{\,2},\, \cdots,\, a_{\,n}\,\right)\,\right\}\,\right\|_{X^{\,\ast}_{d}}.
\end{align*}
Define, \,$L \,:\, X_{F}^{\,\ast} \,\to\, X_{d}^{\,\ast}$\, by
\[L\,(\,T\,) \,=\, \left\{\,T\,\left(\,x_{\,k} \,+\, \lambda\,y_{\,k},\, a_{\,2},\, \cdots,\, a_{\,n}\,\right)\,\right\},\; \;T \,\in\, X_{F}^{\,\ast}.\] 
Then \,$L$\, is a bounded linear operator such that
\begin{align*}
&\left\|\,U\,T \,-\, L\,T\,\right\|_{X^{\,\ast}_{d}} \\
&=\,  \left\|\,\left\{\,T\,\left(\,x_{\,k},\, a_{\,2},\, \cdots,\, a_{\,n}\,\right)\,\right\} \,-\, \left\{\,T\,\left(\,x_{\,k} \,+\, \lambda\,y_{\,k},\, a_{\,2},\, \cdots,\, a_{\,n}\,\right)\,\right\}\,\right\|_{X^{\,\ast}_{d}}\\
&=\, \left\|\,\left\{\,\lambda\,T\,\left(\,y_{\,k},\, a_{\,2},\, \cdots,\, a_{\,n}\,\right)\,\right\}\,\right\|_{X^{\,\ast}_{d}}\\
&\,\leq\, |\,\lambda\,|\,K\,\|\,T\,\|_{X^{\,\ast}_{F}},\; \;T \,\in\, X^{\,\ast}_{F}. 
\end{align*}
This verifies that \,$\|\,U \,-\, L\,\| \,\leq\, |\,\lambda\,|\,K$.\,Now, since \,$S\,U \,=\, I_{F}$, \,$I_{F}$\, is the identity operator on \,$X^{\,\ast}_{F}$, we have 
\begin{align*}
&\left\|\,I_{F} \,-\, S\,L\,\right\| \,=\, \left\|\,S\,U \,-\, S\,L\,\right\| \,\leq\, \|\,S\,\|\,\|\,U \,-\, L\,\| \,<\, 1.  
\end{align*}
Thus \,$S\,L$\, is invertible.\,Take \,$P \,=\, \left(\,S\,L\,\right)^{\,-\, 1}\,S$.\,Then \,$P \,:\, X^{\,\ast}_{d} \,\to\, X_{F}^{\,\ast}$\, is a bounded linear operator such that 
\[P\,\left(\,\left\{\,T\,\left(\,x_{\,k} \,+\, \lambda\,y_{\,k},\, a_{\,2},\, \cdots,\, a_{\,n}\,\right)\,\right\}\,\right) \,=\, T,\; \;T \,\in\, X_{F}^{\,\ast}.\] 
Hence, \,$\left(\,\left\{\,x_{\,k} \,+\, \lambda\,y_{\,k}\,\right\},\, P\,\right)$\, is a retro Banach frame associated to \,$(\,a_{\,2},\, \cdots,\, a_{\,n}\,)$\, for \,$X_{F}^{\,\ast}$\, with respect to \,$X^{\,\ast}_{d}$\, having frame bounds \,$\left(\,\|\,S\,\|^{\,-\, 1} \,-\, |\,\lambda\,|\,K\,\right)$\, and \,$\left(\,B \,+\, |\,\lambda\,|\,K\,\right)$.  
\end{proof}

\begin{theorem}
Let \,$\left(\left\{\,x_{\,k}\,\right\},\, S\,\right)$\, be a retro Banach frame associated to \,$(\,a_{2},\, \cdots,\, a_{n}\,)$\, for \,$X_{F}^{\,\ast}$\, with respect to \,$X^{\,\ast}_{d}$.\,Let \,$\left\{\,y_{\,k}\,\right\} \,\subseteq\, X$\, and \,$\left\{\,\alpha_{\,k}\,\right\} \,\subseteq\, \mathbb{R}$\, be any positively confined sequence such that \,$\left\{\,T\,\left(\,\alpha_{\,k}\,y_{\,k},\, a_{\,2},\, \cdots,\, a_{\,n}\,\right)\,\right\} \,\in\, X^{\,\ast}_{d}$, \,$T \,\in\, X^{\,\ast}_{F}$.\,If \,$V \,:\, X_{F}^{\,\ast} \,\to\, X_{d}^{\,\ast}$\, defined by
\[V\,(\,T\,) \,=\, \left\{\,T\,\left(\,y_{\,k},\, a_{\,2},\, \cdots,\, a_{\,n}\,\right)\,\right\},\; \;T \,\in\, X_{F}^{\,\ast}\]  
such that \,$\|\,V\,\| \,<\, \dfrac{\|\,S\,\|^{\,-\,1}}{\sup\limits_{1 \,\leq\, k \,<\, \infty}\,\alpha_{\,k}}$, then there exists a reconstruction operator \,$P \,:\, X^{\,\ast}_{d} \,\to\, X_{F}^{\,\ast}$\, such that \,$\left(\,\left\{\,x_{\,k} \,+\, \alpha_{\,k}\,y_{\,k}\,\right\},\, P\,\right)$\, is a retro Banach frame associated to \,$(\,a_{\,2},\, \cdots,\, a_{\,n}\,)$\, for \,$X_{F}^{\,\ast}$\, with respect to \,$X^{\,\ast}_{d}$.
\end{theorem}

\begin{proof}
Let \,$U \,:\, X_{F}^{\,\ast} \,\to\, X_{d}^{\,\ast}$\, be a bounded linear operator defined by
\[U\,(\,T\,) \,=\, \left\{\,T\,\left(\,x_{\,k},\, a_{\,2},\, \cdots,\, a_{\,n}\,\right)\,\right\},\; \;T \,\in\, X_{F}^{\,\ast}.\]
It is easy to verify that \,$\left\{\,T\,\left(\,x_{\,k} \,+\, \alpha_{\,k}\,y_{\,k},\, a_{\,2},\, \cdots,\, a_{\,n}\,\right)\,\right\} \,\in\, X^{\,\ast}_{d}$, for all \,$T \,\in\, X^{\,\ast}_{F}$.\,Now, for each \,$T \,\in\, X^{\,\ast}_{F}$, we have
\begin{align*}
&\left\|\,\left\{\,T\,\left(\,x_{\,k} \,+\, \alpha_{\,k}\,y_{\,k},\, a_{\,2},\, \cdots,\, a_{\,n}\,\right)\,\right\}\,\right\|_{X^{\,\ast}_{d}}\\
& \,\leq\, \left\|\,\left\{\,T\,\left(\,x_{\,k},\, a_{\,2},\, \cdots,\, a_{\,n}\,\right)\,\right\}\,\right\|_{X^{\,\ast}_{d}} \,+\, \left\|\,\left\{\,\alpha_{\,k}\,T\,\left(\,y_{\,k},\, a_{\,2},\, \cdots,\, a_{\,n}\,\right)\,\right\}\,\right\|_{X^{\,\ast}_{d}}\\
&\,\leq\, \left\|\,\left\{\,T\,\left(\,x_{\,k},\, a_{\,2},\, \cdots,\, a_{\,n}\,\right)\,\right\}\,\right\|_{X^{\,\ast}_{d}} \,+\, \left(\,\sup\limits_{1 \,\leq\, k \,<\, \infty}\,\alpha_{\,k}\,\right)\,\left\|\,\left\{\,T\,\left(\,y_{\,k},\, a_{\,2},\, \cdots,\, a_{\,n}\,\right)\,\right\}\,\right\|_{X^{\,\ast}_{d}}\\
& \leq\, \left[\,\|\,U\,\| \,+\, \|\,V\,\|\,\left(\,\sup\limits_{1 \,\leq\, k \,<\, \infty}\,\alpha_{\,k}\,\right)\,\right]\,\|\,T\,\|_{X^{\,\ast}_{F}}.
\end{align*}
On the other hand, for each \,$T \,\in\, X^{\,\ast}_{F}$, we have
\begin{align*}
&\left\|\,\left\{\,T\,\left(\,x_{\,k} \,+\, \alpha_{\,k}\,y_{\,k},\, a_{\,2},\, \cdots,\, a_{\,n}\,\right)\,\right\}\,\right\|_{X^{\,\ast}_{d}}\\
&\geq\, \left\|\,\left\{\,T\,\left(\,x_{\,k},\, a_{\,2},\, \cdots,\, a_{\,n}\,\right)\,\right\}\,\right\|_{X^{\,\ast}_{d}} \,-\, \left\|\,\left\{\,\alpha_{\,k}\,T\,\left(\,y_{\,k},\, a_{\,2},\, \cdots,\, a_{\,n}\,\right)\,\right\}\,\right\|_{X^{\,\ast}_{d}}\\
&\geq\, \left[\,\|\,S\,\|^{\,-\, 1} \,-\, \|\,V\,\|\,\left(\,\sup\limits_{1 \,\leq\, k \,<\, \infty}\,\alpha_{\,k}\,\right)\,\right]\,\|\,T\,\|_{X^{\,\ast}_{F}}. 
\end{align*}
Define, \,$L \,:\, X_{F}^{\,\ast} \,\to\, X_{d}^{\,\ast}$\, by
\[L\,(\,T\,) \,=\, \left\{\,T\,\left(\,x_{\,k} \,+\, \alpha_{\,k}\,y_{\,k},\, a_{\,2},\, \cdots,\, a_{\,n}\,\right)\,\right\},\; \;T \,\in\, X_{F}^{\,\ast}.\] 
Then according to proof of the Theorem \ref{th1}, \,$L$\, is a bounded linear operator on \,$X_{F}^{\,\ast}$\, such that \,$\|\,U \,-\, L\,\| \,\leq\, \sup\limits_{1 \,\leq\, k \,<\, \infty}\,\alpha_{\,k}\,\|\,V\,\|$\, and \,$S\,L$\, is invertible.\,Take \,$P \,=\, \left(\,S\,L\,\right)^{\,-\, 1}\,S$.\,Then \,$P \,:\, X^{\,\ast}_{d} \,\to\, X_{F}^{\,\ast}$\, is a bounded linear operator such that 
\[P\,\left(\,\left\{\,T\,\left(\,x_{\,k} \,+\, \alpha_{\,k}\,y_{\,k},\, a_{\,2},\, \cdots,\, a_{\,n}\,\right)\,\right\}\,\right) \,=\, T,\; \;T \,\in\, X_{F}^{\,\ast}.\] 
Hence, \,$\left(\,\left\{\,x_{\,k} \,+\, \alpha_{\,k}\,y_{\,k}\,\right\},\, P\,\right)$\, is a retro Banach frame associated to \,$(\,a_{\,2},\, \cdots,\, a_{\,n}\,)$\, for \,$X_{F}^{\,\ast}$\, with respect to \,$X^{\,\ast}_{d}$. 
\end{proof}

In the next theorem, we establish that retro Banach frame associated to \,$(\,a_{\,2},\, \cdots$, \,$a_{\,n}\,)$\, is stable under perturbation of frame elements by positively confined sequence of scalars.

\begin{theorem}
Let \,$\left(\left\{\,x_{\,k}\,\right\},\, S\,\right)$\, be a retro Banach frame associated to \,$(\,a_{2},\, \cdots,\, a_{n}\,)$\, for \,$X_{F}^{\,\ast}$\, with respect to \,$X^{\,\ast}_{d}$.\,Let \,$\left\{\,y_{\,k}\,\right\} \,\subseteq\, X$\, be such that \,$\left\{\,T\,\left(\,y_{\,k},\, a_{\,2},\, \cdots,\, a_{\,n}\,\right)\,\right\} \,\in\, X^{\,\ast}_{\,d}$, \,$T \,\in\, X^{\,\ast}_{F}$.\,Let \,$R \,:\, X_{d}^{\,\ast} \,\to\, X_{d}^{\,\ast}$\, be a bounded linear operator such that
\[R\,\left(\,\left\{\,T\,\left(\,y_{\,k},\, a_{\,2},\, \cdots,\, a_{\,n}\,\right)\,\right\}\,\right) \,=\, \left\{\,T\,\left(\,x_{\,k},\, a_{\,2},\, \cdots,\, a_{\,n}\,\right)\,\right\},\; \;T \,\in\, X_{F}^{\,\ast}.\]
Suppose \,$\left\{\,\alpha_{\,k}\,\right\}$\, and \,$\left\{\,\beta_{\,k}\,\right\}$\, are two positively confined sequences in \,$\mathbb{R}$.\,If there exist constants \,$\lambda,\,\mu\; \;(\,0 \,\leq\, \lambda,\,\mu \,<\, 1\,)$\, and \,$\gamma$\, such that 
\begin{description}
\item[$(i)$]\,$\gamma \,<\, (\,1 \,-\, \lambda\,)\,\|\,S\,\|^{\,-\,1}\,\left(\,\inf\limits_{1 \,\leq\, k \,<\, \infty}\,\alpha_{\,k}\,\right)$.
\item[$(ii)$]$\left\|\,\left\{\,\alpha_{\,k}\,T\,\left(\,x_{\,k},\, a_{\,2},\, \cdots,\, a_{\,n}\,\right)\,\right\} \,-\, \left\{\,\beta_{\,k}\,T\,\left(\,y_{\,k},\, a_{\,2},\, \cdots,\, a_{\,n}\,\right)\,\right\}\,\right\|_{X^{\,\ast}_{d}}$
\begin{align*}
&\leq\, \lambda\,\left\|\,\left\{\,\alpha_{\,k}\,T\,\left(\,x_{\,k},\, a_{\,2},\, \cdots,\, a_{\,n}\,\right)\,\right\}\,\right\|_{X^{\,\ast}_{d}} \,+\, \mu\,\left\|\,\left\{\,\beta_{\,k}\,T\,\left(\,y_{\,k},\, a_{\,2},\, \cdots,\, a_{\,n}\,\right)\,\right\}\,\right\|_{X^{\,\ast}_{d}} \\
&\hspace{1cm} \,+\, \gamma\,\|\,T\,\|_{X^{\,\ast}_{F}},\; \;T \,\in\, X_{F}^{\,\ast}. 
\end{align*}
\end{description}
Then there exists a reconstruction operator \,$P \,:\, X^{\,\ast}_{d} \,\to\, X_{F}^{\,\ast}$\, such that \,$\left(\,\left\{\,y_{\,k}\,\right\},\, P\,\right)$\, is a retro Banach frame associated to \,$(\,a_{\,2},\, \cdots,\, a_{\,n}\,)$\, for \,$X_{F}^{\,\ast}$\, with respect to \,$X^{\,\ast}_{d}$.
\end{theorem}

\begin{proof}
Let \,$U \,:\, X_{F}^{\,\ast} \,\to\, X_{d}^{\,\ast}$\, be a bounded linear operator defined by
\[U\,(\,T\,) \,=\, \left\{\,T\,\left(\,x_{\,k},\, a_{\,2},\, \cdots,\, a_{\,n}\,\right)\,\right\},\; \;T \,\in\, X_{F}^{\,\ast}.\]
Since the operator \,$S\,U \,:\, X_{F}^{\,\ast} \,\to\, X_{F}^{\,\ast}$\, is an identity operator, for \,$T \,\in\, X_{F}^{\,\ast}$,
\[\|\,T\,\|_{X^{\,\ast}_{F}} \,=\, \|\,S\,U\,(\,T\,)\,\|_{X^{\,\ast}_{F}} \,\leq\, \|\,S\,\|\,\left\|\,\left\{\,T\,\left(\,x_{\,k},\, a_{\,2},\, \cdots,\, a_{\,n}\,\right)\,\right\}\,\right\|_{X^{\,\ast}_{d}}.\]
Now, for each \,$T \,\in\, X_{F}^{\,\ast}$, we have
\begin{align*}
&\left\|\,\left\{\,\beta_{\,k}\,T\,\left(\,y_{\,k},\, a_{\,2},\, \cdots,\, a_{\,n}\,\right)\,\right\}\,\right\|_{X^{\,\ast}_{d}}\\
& \,\leq\, \left\|\,\left\{\,\alpha_{\,k}\,T\,\left(\,x_{\,k},\, a_{\,2},\, \cdots,\, a_{\,n}\,\right)\,\right\}\,\right\|_{X^{\,\ast}_{d}} \,+\\
&\hspace{1cm}\,+\,\left\|\,\left\{\,\alpha_{\,k}\,T\,\left(\,x_{\,k},\, a_{\,2},\, \cdots,\, a_{\,n}\,\right)\,\right\} \,-\, \left\{\,\beta_{\,k}\,T\,\left(\,y_{\,k},\, a_{\,2},\, \cdots,\, a_{\,n}\,\right)\,\right\}\,\right\|_{X^{\,\ast}_{d}} \\
&\leq\, (\,1 \,+\, \lambda\,)\,\left\|\,\left\{\,\alpha_{\,k}\,T\,\left(\,x_{\,k},\, a_{\,2},\, \cdots,\, a_{\,n}\,\right)\,\right\}\,\right\|_{X^{\,\ast}_{d}} \,+\,\\
&\hspace{1cm} \,+\,\mu\,\left\|\,\left\{\,\beta_{\,k}\,T\,\left(\,y_{\,k},\, a_{\,2},\, \cdots,\, a_{\,n}\,\right)\,\right\}\,\right\|_{X^{\,\ast}_{d}} \,+\, \gamma\,\|\,T\,\|_{X^{\,\ast}_{F}}.\\
&\Rightarrow\, (\,1 \,-\, \mu\,)\,\left\|\,\left\{\,\beta_{\,k}\,T\,\left(\,y_{\,k},\, a_{\,2},\, \cdots,\, a_{\,n}\,\right)\,\right\}\,\right\|_{X^{\,\ast}_{d}}\\
& \,\leq\, \left[\,(\,1 \,+\, \lambda\,)\,\|\,U\,\|\,\left(\,\sup\limits_{1 \,\leq\, k \,<\, \infty}\,\alpha_{\,k}\,\right) \,+\, \gamma\,\right]\,\|\,T\,\|_{X^{\,\ast}_{F}}.\\
&\Rightarrow\, (\,1 \,-\, \mu\,)\,\left(\,\inf\limits_{1 \,\leq\, k \,<\, \infty}\,\beta_{\,k}\,\right)\,\left\|\,\left\{\,T\,\left(\,y_{\,k},\, a_{\,2},\, \cdots,\, a_{\,n}\,\right)\,\right\}\,\right\|_{X^{\,\ast}_{d}}\\
& \,\leq\, \left[\,(\,1 \,+\, \lambda\,)\,\|\,U\,\|\,\left(\,\sup\limits_{1 \,\leq\, k \,<\, \infty}\,\alpha_{\,k}\,\right) \,+\, \gamma\,\right]\,\|\,T\,\|_{X^{\,\ast}_{F}}.    
\end{align*}
On the other hand, by condition $(ii)$, we get
\begin{align*}
&(\,1 \,+\, \mu\,)\,\left\|\,\left\{\,\beta_{\,k}\,T\,\left(\,y_{\,k},\, a_{\,2},\, \cdots,\, a_{\,n}\,\right)\,\right\}\,\right\|_{X^{\,\ast}_{d}}\\
&\geq\, (\,1 \,-\, \lambda\,)\,\left\|\,\left\{\,\alpha_{\,k}\,T\,\left(\,x_{\,k},\, a_{\,2},\, \cdots,\, a_{\,n}\,\right)\,\right\}\,\right\|_{X^{\,\ast}_{d}} \,-\, \gamma\,\|\,T\,\|_{X^{\,\ast}_{F}} \\
&\geq\, \left[\,(\,1 \,-\, \lambda\,)\,\|\,S\,\|^{\,-\, 1}\,\left(\,\inf\limits_{1 \,\leq\, k \,<\, \infty}\,\alpha_{\,k}\,\right) \,-\, \gamma\,\right]\,\|\,T\,\|_{X^{\,\ast}_{F}},\; \;T \,\in\, X_{F}^{\,\ast}.
\end{align*}
Therefore, for each \,$T \,\in\, X_{F}^{\,\ast}$, we have
\begin{align*}
&(\,1 \,+\, \mu\,)\,\left(\,\sup\limits_{1 \,\leq\, k \,<\, \infty}\,\beta_{\,k}\,\right)\,\left\|\,\left\{\,T\,\left(\,y_{\,k},\, a_{\,2},\, \cdots,\, a_{\,n}\,\right)\,\right\}\,\right\|_{X^{\,\ast}_{d}}\\
&\geq\, (\,1 \,+\, \mu\,)\,\left\|\,\left\{\,\beta_{\,k}\,T\,\left(\,y_{\,k},\, a_{\,2},\, \cdots,\, a_{\,n}\,\right)\,\right\}\,\right\|_{X^{\,\ast}_{d}}\\
&\geq\, \left[\,(\,1 \,-\, \lambda\,)\,\|\,S\,\|^{\,-\, 1}\,\left(\,\inf\limits_{1 \,\leq\, k \,<\, \infty}\,\alpha_{\,k}\,\right) \,-\, \gamma\,\right]\,\|\,T\,\|_{X^{\,\ast}_{F}}.
\end{align*} 
Thus, for each \,$T \,\in\, X_{F}^{\,\ast}$, we have 
\begin{align*}
&\dfrac{(\,1 \,-\, \lambda\,)\,\|\,S\,\|^{\,-\, 1}\,\left(\,\inf\limits_{1 \,\leq\, k \,<\, \infty}\,\alpha_{\,k}\,\right) \,-\, \gamma}{(\,1 \,+\, \mu\,)\,\left(\,\sup\limits_{1 \,\leq\, k \,<\, \infty}\,\beta_{\,k}\,\right)}\,\|\,T\,\|_{X^{\,\ast}_{F}} \,\leq\, \left\|\,\left\{\,T\,\left(\,y_{\,k},\, a_{\,2},\, \cdots,\, a_{\,n}\,\right)\,\right\}\,\right\|_{X^{\,\ast}_{d}}\\
&\leq\, \dfrac{(\,1 \,+\, \lambda\,)\,\|\,U\,\|\,\left(\,\sup\limits_{1 \,\leq\, k \,<\, \infty}\,\alpha_{\,k}\,\right) \,+\, \gamma}{(\,1 \,-\, \mu\,)\,\left(\,\inf\limits_{1 \,\leq\, k \,<\, \infty}\,\beta_{\,k}\,\right)}\,\|\,T\,\|_{X^{\,\ast}_{F}}. 
\end{align*}
Now, take \,$P \,=\, S\,V$.\,Then \,$P \,:\, X^{\,\ast}_{d} \,\to\, X_{F}^{\,\ast}$\, is a bounded linear operator such that 
\[P\,\left(\,\left\{\,T\,\left(\,y_{\,k},\, a_{\,2},\, \cdots,\, a_{\,n}\,\right)\,\right\}\,\right) \,=\, T,\; \;T \,\in\, X_{F}^{\,\ast}.\] 
Hence, \,$\left(\,\left\{\,y_{\,k}\,\right\},\, P\,\right)$\, is a retro Banach frame associated to \,$(\,a_{\,2},\, \cdots,\, a_{\,n}\,)$\, for \,$X_{F}^{\,\ast}$\, with respect to \,$X^{\,\ast}_{d}$. 
\end{proof}

Next, we give a perturbation result for retro Banach Bessel sequence associated to \,$\left(\,a_{\,2},\, \cdots,\, a_{\,n}\,\right)$.

\begin{theorem}
Let the sequence \,$\left\{\,x_{\,k}\,\right\}$\, be a retro Banach Bessel sequence associated to \,$(\,a_{\,2},\, \cdots,\, a_{\,n}\,)$\, for \,$X_{F}^{\,\ast}$\, with respect to \,$X^{\,\ast}_{d}$.\,Let \,$\left\{\,y_{\,k}\,\right\} \,\subseteq\, X$\, such that \,$\left\{\,T\,\left(\,y_{\,k},\, a_{\,2},\, \cdots,\, a_{\,n}\,\right)\,\right\} \,\in\, X^{\,\ast}_{d}$, \,$T \,\in\, X^{\,\ast}_{F}$.\,Then \,$\left\{\,y_{\,k}\,\right\}$\, is a retro Banach Bessel sequence associated to \,$(\,a_{\,2},\, \cdots,\, a_{\,n}\,)$\, for \,$X_{F}^{\,\ast}$\, with respect to \,$X^{\,\ast}_{d}$, if there exist non-negative constants \,$\alpha,\,\beta\; \;(\,\beta \,<\, 1\,)$\, and \,$\gamma$\, such that for \,$T \,\in\, X_{F}^{\,\ast}$,
\begin{align*}
&\left\|\,\left\{\,T\,\left(\,x_{\,k} \,-\, y_{\,k},\, a_{\,2},\, \cdots,\, a_{\,n}\,\right)\,\right\}\,\right\|_{X^{\,\ast}_{d}} \,\leq\, \alpha\, \left\|\,\left\{\,T\,\left(\,x_{\,k},\, a_{\,2},\, \cdots,\, a_{\,n}\,\right)\,\right\}\,\right\|_{X^{\,\ast}_{d}}\\
&\hspace{1cm} \,+\, \beta\,\left\|\,\left\{\,T\,\left(\,y_{\,k},\, a_{\,2},\, \cdots,\, a_{\,n}\,\right)\,\right\}\,\right\|_{X^{\,\ast}_{d}} \,+\, \gamma\,\|\,T\,\|_{X^{\,\ast}_{F}}.
\end{align*}
\end{theorem} 

\begin{proof}
Since \,$\left\{\,x_{\,k}\,\right\}$\, be a retro Banach Bessel sequence associated to \,$(\,a_{\,2},\, \cdots,\, a_{\,n}\,)$,
\[\left\|\,\left\{\,T\,\left(\,x_{\,k},\, a_{\,2},\, \cdots,\, a_{\,n}\,\right)\,\right\}\,\right\|_{X^{\,\ast}_{d}} \,\leq\, B\,\|\,T\,\|_{X^{\,\ast}_{F}}\; \;\forall\; \,T \,\in\, X^{\,\ast}_{F}. \]
Thus, for each \,$T \,\in\, X^{\,\ast}_{F}$, we have
\begin{align*}
&\left\|\,\left\{\,T\,\left(\,y_{\,k},\, a_{\,2},\, \cdots,\, a_{\,n}\,\right)\,\right\}\,\right\|_{X^{\,\ast}_{d}}\\
&=\, \left\|\,\left\{\,T\,\left(\,x_{\,k},\, a_{\,2},\, \cdots,\, a_{\,n}\,\right)\,\right\} \,-\, \left\{\,T\,\left(\,x_{\,k} \,-\, y_{\,k},\, a_{\,2},\, \cdots,\, a_{\,n}\,\right)\,\right\}\,\right\|_{X^{\,\ast}_{d}}\\
&\leq\, (\,1 \,+\, \alpha\,)\,\left\|\,\left\{\,T\,\left(\,x_{\,k},\, a_{\,2},\, \cdots,\, a_{\,n}\,\right)\,\right\}\,\right\|_{X^{\,\ast}_{d}} \,+\, \beta\,\left\|\,\left\{\,T\,\left(\,y_{\,k},\, a_{\,2},\, \cdots,\, a_{\,n}\,\right)\,\right\}\,\right\|_{X^{\,\ast}_{d}}\\
&\hspace{1cm}\,+\, \gamma\,\|\,T\,\|_{X^{\,\ast}_{F}}\\
&\leq\, \dfrac{(\,1 \,+\, \alpha\,)\,B \,+\, \gamma}{1 \,-\, \beta}\,\|\,T\,\|_{X^{\,\ast}_{F}}.
\end{align*} 
Hence, \,$\left\{\,y_{\,k}\,\right\}$\, is a retro Banach Bessel sequence associated to \,$(\,a_{\,2},\, \cdots,\, a_{\,n}\,)$\, for \,$X_{F}^{\,\ast}$\, with respect to \,$X^{\,\ast}_{d}$.  
\end{proof}

The following two theorems gives a condition under which the finite sum of retro Banach frame associated to \,$(\,a_{\,2},\, \cdots,\, a_{\,n}\,)$\, is again a retro Banach frame associated to \,$(\,a_{\,2},\, \cdots,\, a_{\,n}\,)$.

\begin{theorem}
For \,$i \,\in\, \left\{\,1,\,2,\,\cdots,\,m\,\right\}$, let \,$\left(\,\left\{\,x_{\,k,\,i}\,\right\},\, S_{\,i}\,\right)$\, be retro Banach frames associated to \,$(\,a_{\,2},\, \cdots$, \,$a_{\,n}\,)$\, for \,$X_{F}^{\,\ast}$\, with respect to \,$X^{\,\ast}_{d}$.\,Then there exists a reconstruction operator \,$P \,:\, X^{\,\ast}_{d} \,\to\, X_{F}^{\,\ast}$\, such that \,$\left(\,\left\{\,\sum\limits_{i \,=\, 1}^{m}\,x_{\,k,\,i}\,\right\},\, P\,\right)$\, is a tight retro Banach frame associated to \,$(\,a_{\,2},\, \cdots,\, a_{\,n}\,)$\, for \,$X_{F}^{\,\ast}$\, with respect to \,$X^{\,\ast}_{d}$, provided
\[\left\|\,\left\{\,T\,\left(\,x_{\,k,\,j},\, a_{\,2},\, \cdots,\, a_{\,n}\,\right)\,\right\}\,\right\|_{X^{\,\ast}_{d}} \,\leq\, \left\|\,\left\{\,T\,\left(\,\sum\limits_{i \,=\, 1}^{m}\,x_{\,k,\,i},\, a_{\,2},\, \cdots,\, a_{\,n}\,\right)\,\right\}\,\right\|_{X^{\,\ast}_{d}},\] 
for \,$T \,\in\, X_{F}^{\,\ast}$, for some \,$j \,\in\, \left\{\,1,\, 2,\, \cdots,\, m\,\right\}$. 
\end{theorem}

\begin{proof}
For \,$T \,\in\, X_{F}^{\,\ast}$, we have
\begin{align*}
\|\,T\,\|_{X^{\,\ast}_{F}} &\,=\, \left\|\,S_{\,i}\,\left(\,\left\{\,T\,\left(\,x_{\,k,\,j},\, a_{\,2},\, \cdots,\, a_{\,n}\,\right)\,\right\}\,\right)\,\right\|_{X^{\,\ast}_{F}}\\
&\leq\, \left\|\,S_{\,i}\,\right\|\,\left\|\,\left\{\,T\,\left(\,x_{\,k,\,j},\, a_{\,2},\, \cdots,\, a_{\,n}\,\right)\,\right\}\,\right\|_{X^{\,\ast}_{d}}\\
&\leq\,\left\|\,S_{\,i}\,\right\|\,\left\|\,\left\{\,T\,\left(\,\sum\limits_{i \,=\, 1}^{m}\,x_{\,k,\,i},\, a_{\,2},\, \cdots,\, a_{\,n}\,\right)\,\right\}\,\right\|_{X^{\,\ast}_{d}}. 
\end{align*} 
Thus, \,$\left\{\,T\,\left(\,\sum\limits_{i \,=\, 1}^{m}\,x_{\,k,\,i},\, a_{\,2},\, \cdots,\, a_{\,n}\,\right)\,\right\}$\, is total over \,$X_{F}^{\,\ast}$.\,Therefore, by Remark 7.1 in \cite{I}, there exists an associated Banach space 
\[X_{d_{\,1}} \,=\, \left\{\,\left\{\,T\,\left(\,\sum\limits_{i \,=\, 1}^{m}\,x_{\,k,\,i},\, a_{\,2},\, \cdots,\, a_{\,n}\,\right)\,\right\} \,:\, T \,\in\, X_{F}^{\,\ast}\,\right\}\]
equipped with the norm 
\[\left\|\,\left\{\,T\,\left(\,\sum\limits_{i \,=\, 1}^{m}\,x_{\,k,\,i},\, a_{\,2},\, \cdots,\, a_{\,n}\,\right)\,\right\}\,\right\|_{X_{d_{\,1}}} \,=\, \|\,T\,\|_{X^{\,\ast}_{F}},\; \;T \,\in\, X_{F}^{\,\ast}\]
and a bounded linear operator \,$P \,:\, X^{\,\ast}_{\,d} \,\to\, X_{F}^{\,\ast}$\, defined by
\[P\,\left(\,\left\{\,T\,\left(\,\sum\limits_{i \,=\, 1}^{m}\,x_{\,k,\,i},\, a_{\,2},\, \cdots,\, a_{\,n}\,\right)\,\right\}\,\right) \,=\, T,\; \;T \,\in\, X_{F}^{\,\ast}\]
such that \,$\left(\,\left\{\,\sum\limits_{i \,=\, 1}^{m}\,x_{\,k,\,i}\,\right\},\, P\,\right)$\, is a tight retro Banach frame associated to \,$(\,a_{\,2},\, \cdots,$ \,$a_{\,n}\,)$\, for \,$X_{F}^{\,\ast}$\, with respect to \,$X^{\,\ast}_{d}$. 
\end{proof}

\begin{theorem}
Let \,$\left(\,\left\{\,x_{\,k,\,i}\,\right\},\, S_{\,i}\,\right),\; \;i \,\in\, E_{\,m} \,=\, \left\{\,1,\,2,\,\cdots,\,m\,\right\}$\, be retro Banach frames associated to \,$(\,a_{\,2},\, \cdots$, \,$a_{\,n}\,)$\, for \,$X_{F}^{\,\ast}$\, with respect to \,$X^{\,\ast}_{d}$.\,Let \,$\left\{\,y_{\,k,\,i}\,\right\} \,\subseteq\, X$\, be such that \,$\left\{\,T\,\left(\,y_{\,k,\,i},\, a_{2},\, \cdots,\, a_{n}\,\right)\,\right\} \,\in\, X^{\,\ast}_{d}$, \,$T \,\in\, X^{\,\ast}_{F}$.\,Suppose \,$R \,:\, X_{d}^{\,\ast} \,\to\, X_{d}^{\,\ast}$\, be a bounded linear operator such that
\[R\,\left(\,\left\{\,T\,\left(\,\sum\limits_{i \,=\, 1}^{\,m}\,y_{\,k,\,i},\, a_{\,2},\, \cdots,\, a_{\,n}\,\right)\,\right\}\,\right) \,=\, \left\{\,T\,\left(\,x_{\,k,\,p},\, a_{\,2},\, \cdots,\, a_{\,n}\,\right)\,\right\},\; \;T \,\in\, X_{F}^{\,\ast},\]
for some \,$p \,\in\, E_{\,m}$\, and for each \,$i \,=\, 1,\,2,\,\cdots,\,m$, let \,$U_{\,i} \,:\, X_{F}^{\,\ast} \,\to\, X_{d}^{\,\ast}$\, be an operator defined by
\[U_{\,i}\,(\,T\,) \,=\, \left\{\,T\,\left(\,x_{k,\,i},\, a_{\,2},\, \cdots,\, a_{\,n}\,\right)\,\right\},\; \;T \,\in\, X_{F}^{\,\ast}.\] 
If there exist constants \,$\alpha,\,\beta \,>\, 0$\, such that
\begin{description}
\item[$(i)$]\hspace{.5cm}$\alpha\,\sum\limits_{i \,\in\, E_{\,m}}\,\left\|\,U_{\,i}\,\right\| \,+\, m\,\beta \,<\, \left\|\,S_{j}\,\right\|^{\,-\, 1} \,-\, \sum\limits_{i \,\in\, E_{\,m},\, i \,\neq\, j}\,\left\|\,U_{\,i}\,\right\|$, for some \,$j \,\in\, E_{\,m}$.
\item[$(ii)$]\hspace{.5cm}$\left\|\,\left\{\,T\,\left(\,x_{\,k,\,i} \,-\, y_{\,k,\,i},\, a_{\,2},\, \cdots,\, a_{\,n}\,\right)\,\right\}\,\right\|_{X^{\,\ast}_{d}}$
\begin{align*}
& \,\leq\, \alpha\, \left\|\,\left\{\,T\,\left(\,x_{\,k,\,i},\, a_{\,2},\, \cdots,\, a_{\,n}\,\right)\,\right\}\,\right\|_{X^{\,\ast}_{d}} \,+\, \beta\,\|\,T\,\|_{X^{\,\ast}_{F}},\; \;T \,\in\, X_{F}^{\,\ast},\; \;i \,\in\, E_{\,m},
\end{align*}
\end{description}
then there exists a bounded linear operator \,$P \,:\, X^{\,\ast}_{d} \,\to\, X_{F}^{\,\ast}$\, such that the family \,$\left(\,\left\{\,\sum\limits_{i \,\in\, E_{\,m}}\,y_{\,k,\,i}\,\right\},\, P\,\right)$\, is a retro Banach frame associated to \,$(\,a_{\,2},\, \cdots,\, a_{\,n}\,)$\, for \,$X_{F}^{\,\ast}$\, with respect to \,$X^{\,\ast}_{d}$.
\end{theorem}

\begin{proof}
For each \,$i \,\in\, E_{\,m}$, \,$S_{\,i}\,U_{\,i}$\, is an identity operator on \,$X_{F}^{\,\ast}$.\,Therefore, for each \,$T \,\in\, X_{F}^{\,\ast}$, we have
\begin{equation}\label{eq1.41}
\|\,T\,\|_{X^{\,\ast}_{F}} \,=\, \left\|\,S_{\,i}\,U_{\,i}\,T\,\right\|_{X^{\,\ast}_{F}} \,\leq\, \|\,S_{\,i}\,\|\,\left\|\,\left\{\,T\,\left(\,x_{\,k,\,i},\, a_{\,2},\, \cdots,\, a_{\,n}\,\right)\,\right\}\,\right\|_{X^{\,\ast}_{d}}.
\end{equation}
Also, for \,$T \,\in\, X_{F}^{\,\ast}$, we have
\begin{align}\label{eq1.51}
\left\|\,\left\{\,T\,\left(\,\sum\limits_{i \,\in\, E_{\,m}}\,x_{\,k,\,i},\, a_{\,2},\, \cdots,\, a_{\,n}\,\right)\,\right\}\,\right\|_{X^{\,\ast}_{d}} &\,=\, \left\|\,\sum\limits_{i \,\in\, E_{\,m}}\,\left\{\,T\,\left(\,x_{\,k,\,i},\, a_{\,2},\, \cdots,\, a_{\,n}\,\right)\,\right\}\,\right\|_{X^{\,\ast}_{d}}\nonumber\\
&\leq\, \sum\limits_{i \,\in\, E_{\,m}}\,\left\|\,U_{\,i}\,\right\|\,\|\,T\,\|_{X^{\,\ast}_{F}}.  
\end{align}
Now, for each \,$T \,\in\, X_{F}^{\,\ast}$, we have
\begin{align*}
&\left\|\,\left\{\,T\,\left(\,\sum\limits_{i \,\in\, E_{\,m}}\,y_{\,k,\,i},\, a_{\,2},\, \cdots,\, a_{\,n}\,\right)\,\right\}\,\right\|_{X^{\,\ast}_{d}}\\
&=\,\left\|\,\sum\limits_{i \,\in\, E_{\,m}}\,\left\{\,T\,\left(\,x_{\,k,\,i},\, a_{\,2},\, \cdots,\, a_{\,n}\,\right) \,-\, T\,\left(\,x_{\,k,\,i} \,-\, y_{\,k,\,i},\, a_{\,2},\, \cdots,\, a_{\,n}\,\right) \,\right\}\,\right\|_{X^{\,\ast}_{d}}\\
&\geq\, \left\|\,\sum\limits_{i \,\in\, E_{\,m}}\,\left\{\,T\,\left(\,x_{\,k,\,i},\, a_{\,2},\, \cdots,\, a_{\,n}\,\right)\,\right\}\,\right\|_{X^{\,\ast}_{\,d}} \,-\, \left\|\,\sum\limits_{i \,\in\, E_{\,m}}\,\left\{\,T\,\left(\,x_{\,k,\,i} \,-\, y_{\,k,\,i},\, a_{\,2},\, \cdots,\, a_{\,n}\,\right) \,\right\}\,\right\|_{X^{\,\ast}_{d}}\\
&\geq\, \left\|\,\left\{\,T\,\left(\,x_{\,k,\,j},\, a_{\,2},\, \cdots,\, a_{\,n}\,\right)\,\right\} \,+\, \sum\limits_{i \,\in\, E_{\,m},\,i \,\neq\, j}\,\left\{\,T\,\left(\,x_{\,k,\,i},\, a_{\,2},\, \cdots,\, a_{\,n}\,\right)\,\right\}\,\right\|_{X^{\,\ast}_{d}} \,-\\
&\hspace{1cm} \,-\, \left\|\,\sum\limits_{i \,\in\, E_{\,m}}\,\left\{\,T\,\left(\,x_{\,k,\,i} \,-\, y_{\,k,\,i},\, a_{\,2},\, \cdots,\, a_{\,n}\,\right) \,\right\}\,\right\|_{X^{\,\ast}_{d}}\\ 
&\geq\, \left\|\,S_{j}\,\right\|^{\,-\, 1} \,-\, \left[\,\alpha\,\sum\limits_{i \,\in\, E_{\,m}}\,\left\|\,U_{\,i}\,\right\| \,+\, \sum\limits_{i \,\in\, E_{\,m},\, i \,\neq\, j}\,\left\|\,U_{\,i}\,\right\| \,+\, m\,\beta\,\right]\,\|\,T\,\|_{X^{\,\ast}_{F}}\; \;[\;\text{by (\ref{eq1.41}) and (\ref{eq1.51})}\;].    
\end{align*}
On the other hand, using (\ref{eq1.51}), for each \,$T \,\in\, X_{F}^{\,\ast}$, we get
\[\left\|\left\{\,T\left(\,\sum\limits_{i \,\in\, E_{\,m}}y_{\,k,\,i},\, a_{2},\, \cdots,\, a_{n}\,\right)\right\}\right\|_{X^{\,\ast}_{d}} \,\leq\, \left(\,(\,1 \,+\, \alpha\,)\,\sum\limits_{i \,\in\, E_{\,m}}\,\left\|\,U_{\,i}\,\right\| \,+\, \,\beta\,\right)\,\|\,T\,\|_{X^{\,\ast}_{F}}. \]
Now, we take \,$P \,=\, S_{\,p}\,R$, where \,$p$\, is fixed.\,Then \,$P \,:\, X^{\,\ast}_{d} \,\to\, X_{F}^{\,\ast}$\, is a bounded linear operator such that
\[P\,\left(\,\left\{\,T\,\left(\,\sum\limits_{i \,\in\, E_{\,m}}\,y_{\,k,\,i},\, a_{\,2},\, \cdots,\, a_{\,n}\,\right)\,\right\}\,\right) \,=\, T,\; \;T \,\in\, X_{F}^{\,\ast}.\]
Hence, \,$\left(\,\left\{\,\sum\limits_{i \,\in\, E_{\,m}}\,y_{\,k,\,i}\,\right\},\, P\,\right)$\, is a retro Banach frame associated to \,$(\,a_{\,2},\, \cdots,\, a_{\,n}\,)$\, for \,$X_{F}^{\,\ast}$\, with respect to \,$X^{\,\ast}_{d}$. 
\end{proof}

\end{document}